\documentclass[12pt]{article}
\usepackage{amsmath,amsthm,amssymb, amstext, amsopn, amsxtra}
\usepackage[all]{xy}
\usepackage{latexsym}
\usepackage{amsfonts}
\begin{document}
\newcommand{\nt}{\noindent}
\newcommand{\bs}{\bigskip}
\newcommand{\ms}{\medskip}
\newcommand{\mk}{\medskip}
\newcommand{\sk}{\smallskip}
\newcommand{\ep}{\varepsilon}
\newcommand{\m}{{\mathfrak m}}
\newcommand{\p}{{\mathfrak p}}
\newcommand{\dd}{\delta}
\newcommand{\A}{{\mathbb A}}
\newcommand{\R}{{\mathcal R}}
\newcommand{\M}{{\mathbb M}}
\newcommand{\Rstar}{$\,(R,\,{}^*\,)\,$}
\newcommand{\starring}{$\,{\ast}\,$-ring}
\newcommand{\starrings}{$\,{\ast}\,$-rings}
\newcommand{\ap}{\alpha^{\,\prime}}
\newcommand{\bp}{\beta^{\,\prime}}
\newcommand{\isom}{\cong}




\bs

\sk
\begin{center}
\Large{\bf Ring Elements of Stable Range One}

\bigskip
\large{Dinesh Khurana and T.\,Y.~Lam}
\end{center}

\mk
\begin{abstract}
\begin{small}
\nt
A ring element $\,a\in R\,$ is said to be of {\it right stable range
one\/} if, for any $\,t\in R$, $\,aR+tR=R\,$ implies that $\,a+t\,b\,$
is a unit in $\,R\,$ for some $\,b\in R$.  Similarly, $\,a\in R\,$
is said to be of {\it left stable range one\/} if $\,R\,a+R\,t=R\,$
implies that $\,a+b't\,$ is a unit in $\,R\,$ for some $\,b'\in R$.
In the last two decades, it has often been speculated that these two
notions are actually the same for any $\,a\in R$. In \S3 of this paper,
we will prove that this is indeed the case.  The key to the proof
of this new symmetry result is a certain ``Super Jacobson's Lemma'',
which generalizes Jacobson's classical lemma stating that, for any
$\,a,b\in R$, $\,1-ab\,$ is a unit in $\,R\,$ iff so is $\,1-ba$.  Our
proof for the symmetry result above has led to a new generalization of
a classical determinantal identity of Sylvester, which will be published
separately in [KL$_3$].  In \S\S4-5, a detailed study is offered for
stable range one ring elements that are unit-regular or nilpotent,
while \S6 examines the behavior of stable range one elements via their
classical Peirce decompositions.  The paper ends with a more concrete
\S7 on integral matrices of stable range one, followed by a final \S8
with a few open questions.
\end{small}
\end{abstract}

\nt
{\bf \S1. \ Introduction}

\bs\indent
The concept and the basic theory of the stable range of unital rings
were first introduced by H.~Bass [Ba] in 1964 in his original study
of the stability properties of the general linear groups in algebraic
K-theory. Over the years, many papers have been written on the stable
range of rings, while the notion of {\it rings of stable range one\/}
has often stood out as a particularly important special
case. In 2005, inspired by Bass's ground-breaking paper, the authors
introduced in [KL$_2$:~Def.~(3.1)] an {\it element-wise\/} notion of
stable range one, as follows.

\bs\nt
{\bf Definition 1.1.} We say that a ring element $\,a\in R\,$ has
{\it right stable range one\/} (written $\,{\rm sr}_R(a)=1$, or
just $\,{\rm sr}\,(a)=1\,$ if $\,R\,$ is understood) if, for any
$\,t\in R$, $\,aR+tR=R\,$ implies that $\,a+t\,b\in {\rm U}(R)\,$
(the group of units in $\,R\,$) for some $\,b\in R$.  Similarly,
we say that $\,a\in R\,$ has {\it left stable range one\/} (written
$\,{\rm sr}\,'(a)={\rm sr}\,'_{R}(a)=1$\/) if, for any $\,t\in R$,
$\,R\,a+R\,t=R\,$ implies that $\,a+b'\,t \in {\rm U}(R)\,$ for
some $\,b'\in R$. Finally, the ring $\,R\,$ itself is said to be
of stable range one if $\,{\rm sr}\,(a)=1\,$ for all $\,a\in R$.


\bs
From the abstract definitions above, it is of course not clear at all
whether, for any element $\,a\in R$, $\,{\rm sr}\,(a)=1\,$ is equivalent
to $\,{\rm sr}\,'(a)=1\,$ (over any ring $\,R$). Quite pleasantly, this
long-suspected statement did turn out to be true: we will give a nontrivial
proof for it in \S3 after covering some general preparatory material
in \S2. After this, of course, the distinction between the two conditions
$\,{\rm sr}\,(a)=1\,$ and $\,{\rm sr}\,'(a)=1\,$ will become unnecessary.
On the other hand, if $\,a\in S\,$ where $\,S\,$ is a subring of $\,R$,
there does not seem to be much correlation between the condition
$\,{\rm sr}_S(a)=1\,$ and the condition $\,{\rm sr}_R(a)=1$. Nevertheless,
if $\,f:R\to S\,$ is a {\it surjective\/} ring homomorphism, then
$\,{\rm sr}_R(a)=1\,$ does imply that $\,{\rm sr}_S(f(a))=1$; see
Theorem 2.4(D) below.

\bs
This paper is intended to be an introduction to the notion of
``element-wise stable range one'' essentially from first principles.
In particular, a couple of results in the literature that are more
or less known to experts (e.g.~the Product Theorem 2.8) are re-proved in
\S2 and \S5 for the reader's convenience. For a sampling of old and new
papers in the literature written on the theme of stable range one for
rings and for ring elements, we can cite for instance [EO], [Ka], [Ev],
[Va$_1$, Va$_2$], [Wa], [GM], [CY$_1$, CY$_2$], [Yu$_1$, Yu$_2$], [Can],
[Ar], [Ch$_3$, Ch$_4$], [We], [CN], [Khu], [KhM], [AO$_1$, AO$_2$],
[HN], [CLM], [DGK], and [CP$_1$, CP$_2$]. As we have mentioned in the
previous paragraph, the first main result of our paper is the following.

\bs\nt
{\bf Theorem 1.2.} {\it For any ring element $\,a\in R$,
$\,{\rm sr}\,(a)=1\,$ iff $\,{\rm sr}\,'(a)=1$.}

\bs
The key tool to be used for proving this important symmetry result
is {\bf Super Jacobson's Lemma 3.2}, which in its most primitive
form states that, for any three elements $\,a,b,x\,$ in a ring
$\,R\,$ with unit group $\,{\rm U}(R)\,$:
$$
a+b-axb\in {\rm U}(R)\;\,\Longleftrightarrow \;\,a+b-bxa\in {\rm U}(R).
\leqno (1.3)
$$
In fact, in the full version of Super Jacobson's Lemma 3.2, the
equivalence in (1.3) above is shown to hold if $\,{\rm U}(R)\,$
is replaced by the sets of regular or unit-regular elements in any
ring $\,R$. (Later in Theorem 6.10, Super Jacobson's Lemma will also
be proved if $\,{\rm U}(R)\,$ in (1.3) is replaced by the set of
stable range one elements in $\,R$.)  The detailed proofs of these
results will be given in \S3, after some basic facts on ``element-wise
stable range one'' (such as multiplicative closure) are recalled
or proved in \S2.

\bs
In \S4, we study in detail the relationship between stable range
one elements and unit-regular elements in any ring $\,R$.  While
it is essentially well known (say from [Go: Proposition 4.12] or
[KL$_2$: Theorem 3.5]) that unit-regular elements in $\,R\,$ are
exactly the (von Neumann) regular elements in $\,R\,$ with stable
range one, a few more characterizations of unit-regular elements
in terms of the stable range one notion are given in Theorem 4.3.
As an application, it is shown for instance in Corollary 4.7 that
a ring $\,R\,$ is an IC ring (``ring with internal cancellation'')
in the sense of [KL$_2$] iff every regular element in $\,R\,$ has
stable range one, or equivalently, is a product of unit-regular
elements.

\bs
Instead of working with regular and unit-regular elements, \S5 is
devoted to the study of nilpotent elements in a ring $\,R$.  If
$\,a\in R\,$ is a {\it strongly nilpotent element\/} (in the sense
of Levitzki [Lv]), or a {\it central quasi-nilpotent element\/}
(in the sense of Harte [Ha]), it is shown in Theorem 5.1 that we must
have $\,{\rm sr}\,(a)=1$.  In the case of an exchange ring $\,R\,$
(in the sense of Warfield), we show in Theorem 5.2 that every regular
nilpotent element $\,a\in R\,$ is unit-regular, with the property
that $\,{\rm sr}\,(a^i)=1\,$ for all $\,i\geq 1$.  For any pair
of orthogonal idempotents $\,e,f\,$ in any ring $\,R$, Theorem 5.5
shows that $\,{\rm sr}\,(erf)=1\,$ for all $\,r\in R$, and also that
$\,{\rm sr}\,(erfse)=1\,$ for all $\,r,s\in R$.  For some concrete
examples, it is shown as a part of Theorem 5.8 that, if $\,R\,$
is a matrix ring $\,{\mathbb M}_n(S)\,$ where $\,n\geq 2$, then
$\,{\rm sr}\,\bigl(aE_{kj}\bigr)=1\,$ for all $\,a\in S\,$ and all
$\,k,\,j\,$, where $\,\{E_{kj}\}\,$ are the matrix units in $\,R$.

\bs
As is to be expected, the element-wise notion of stable range one
behaves well with respect to the passage from a ring $\,R\,$ to its
Peirce corner rings $\,eRe\,$ for any idempotent $\,e\in R$. In \S6,
we prove a Suspension Theorem 6.2 to the effect that, if $\,f\!:=1-e\,$
and $\,a\in eR\,e$, $\,p\in fR\,e$, then $\,{\rm sr}_{eR\,e}(a)=1\,$
iff $\,{\rm sr}_R(a+p+f)=1$, in which case we have automatically
$\,{\rm sr}_R(a)=1$.  In the case where $\,R\,$ is a matrix ring
$\,{\mathbb M}_n(S)\,$ and $\,A=(a_{ij})\in R\,$ is a lower triangular
matrix with $\,{\rm sr}_S(a_{ii})=1\,$ for all $\,i\,$, Theorem 6.8
shows that $\,{\rm sr}_R(A)=1$, while Theorem 7.1 shows that if $\,A\,$
is a diagonal matrix with at least one diagonal entry $\,a_{ii}=0$,
then again $\,{\rm sr}_R(A)=1$.  Finally, if $\,S\,$ is assumed to
be a commutative elementary divisor domain, Theorem 7.2 shows that
$\,{\rm sr}_R(B)=1\,$ for any matrix $\,B\in {\mathbb M}_n(S)\,$
such that $\,{\rm det}\,(B)\,$ is either zero or a unit in $\,S$.
After \S7, the paper concludes with a final section \S8 which lists
several open problems on stable range one elements in rings.

\bs
The terminology and notations introduced so far in this Introduction
will be used freely throughout the paper. In addition, for any ring
$\,R$, we shall use the standard notations $\,{\rm idem}\,(R)$,
$\,{\rm nil}\,(R)$, $\,{\rm reg}\,(R)$, $\,{\rm ureg}\,(R)\,$ and
$\,{\rm sreg}\,(R)\,$ respectively for the sets of idempotents,
nilpotent elements, (von Neumann) regular elements, unit-regular
elements, and strongly regular elements in any ring $\,R$. Also,
as we have mentioned before, we will use $\,E_{ij}\,$ to denote the
matrix units in any matrix ring $\,{\mathbb M}_n(S)\,$ over a given
ring $\,S$. Other standard terminology and conventions in ring theory
follow mainly those in [Go], [La$_1$] and [La$_2$].

\mk

\bs\nt
{\bf \S2. \ Basic Facts On Element-Wise Stable Range One}

\bs
Throughout this section (and all subsequent sections), the term
``stable range one'' shall mean ``right stable range one'', and
the notation $\,{\rm sr}\,(a)={\rm sr}_R(a)=1\,$ will be taken
to mean that the element $\,a\,$ has right stable range one in
$\,R$.  We begin by mentioning some of the most basic examples
of such elements.

\bs\nt
{\bf Example 2.1.} To begin with, if $\,a\in {\rm U}(R)$, clearly
$\,{\rm sr}\,(a)=1\,$ since we can pick $\,b=0\,$ in Definition 1.1.
(A more general statement on 1-sided invertible elements will
be given in Theorem 2.4(C) below.)  Next, for any $\,e\in
{\rm idem}\,(R)\,$ with complementary idempotent $\,f\!:=1-e$,
we have $\,{\rm sr}\,(e)=1$. Indeed, if $\,e\,x+ty=1\,$ for some
elements $\,t,x,y\in R$, then $\,exf+tyf=f$. Adding $\,e\,$ gives
$\,e+tyf=1-exf\in {\rm U}(R)\,$ (with inverse $\,1+exf$), so we
have checked that $\,{\rm sr}\,(e)=1$.


\bs
The following result offers a number of alternative characterizations
for elements of stable range one in any ring $\,R$.

\bs\nt
{\bf Theorem 2.2.} {\it For any element $\,a\in R$, the
following statements are equivalent\/}:

\mk\nt
(1) $\,{\rm sr}\,(a)=1$.\\
(2) {\it $\,aR+K=R\,$ for a right ideal $\,K\subseteq R
\,\Rightarrow\,a+k\in {\rm U}(R)\,$ for some $\,k\in K$.}\\
(3) {\it For any $\,x\in R$, there exists $\,b\in R\,$ such that
$\,a+b-axb \in {\rm U}(R)$.}\\
(4) {\it For any $\,x,c\in R$, $\,ax+c\in {\rm U}(R)\,\Rightarrow\,
a+c\,b\in {\rm U}(R)\,$ for some $\,b\in R$.}\\
(5) {\it For any $\,x,c\in R$, $\,ax+c=1\,\Rightarrow\,
a+c\,b\in {\rm U}(R)\,$ for some $\,b\in R$.}

\bs\nt
{\bf Proof.} $(1)\Rightarrow (2)\Rightarrow (3)\,$ are trivial
implications.

\mk\nt
$(3)\Rightarrow (4)$. Assume (3), and consider any $\,x,c\in R\,$
such that $\,u\!:=ax+c\in {\rm U}(R)$. Then $\,axu^{-1}+cu^{-1}=1$.
Applying (3), we have $\,a+(1-axu^{-1})\,b_0\in {\rm U}(R)\,$ for
some $\,b_0\in R$. Thus, $\,a+c\,(u^{-1}b_0)\in {\rm U}(R)$, which
proves (4).

\mk\nt
$(4)\Rightarrow (5)\,$ is a tautology.

\mk\nt
$(5)\Rightarrow (1)$. Assume (5), and consider any $\,t\in R\,$
such that $\,aR+tR=R$. Fixing an equation $\,ax+ty=1\,$ and
applying (5) with $\,c=ty$, we see that $\,a+(ty)\,b \in
{\rm U}(R)\,$ for some $\,b\in R$. According to Definition 1.1,
this checks that $\,{\rm sr}\,(a)=1$. \qed

\bs
In the special case where the right ideal $\,aR\,$ is subject to
some ``suitability condition'' in the sense of Nicholson ([Ni$_1$],
[Ni$_2$]), Theorem 2.2 can be transformed in some way into the
corollary below, which is directly inspired by the work of
Camillo-Yu (in [CY$_1$:~Lemma 2] and [Yu$_1$:~Theorem 9]) on the
criterion for an exchange ring to have stable range one.

\bs\nt
{\bf Corollary 2.3.} {\it Let $\,a\in R\,$ be such that, whenever
$\,aR+tR=R\,$ for some $\,t\in R$, there exists an idempotent
$\,f\in aR\,$ such that $\,1-f\in tR$.} ({\it This would be the case,
for instance, if all elements in $\,aR\,$ are suitable in $\,R\,$};
{\it see\/} [DKN:~Example 3.2].) {\it Then the following three
statements are equivalent\,}:


\mk\nt
(1) $\,{\rm sr}\,(a)=1$.\\
(6) {\it For any $\,x\in R\,$ such that $\,ax=(ax)^2$, there
exists $\,b\in R$ such that $\,a+(1-ax)\,b\in {\rm U}(R)$.}\\
(7) {\it For any $\,x\in R\,$ and $\,e\in {\rm idem}\,(R)$,
$\,ax+e=1\,\Rightarrow\,a+e\,b\in {\rm U}(R)\,$ for some $\,b\in R$.}

\bs\nt
{\bf Proof.} Referring back to condition (5) in Theorem 2.2,
it  suffices to check that $(5)\Rightarrow (7)\Rightarrow (6)
\Rightarrow (1)$. Here, the first two implications are obvious.
To prove the last implication, assume that $(6)$ holds. To check (1),
we consider any equation $\,aR+tR=R$. By our hypothesis on $\,a$,
there exists an idempotent $\,f=ax\,$ for some $\,x\in R\,$ such
that $\,e\!:=1-f=ty\,$ for some $\,y\in R$. Since $\,ax+e=f+e=1$,
$(6)$ guarantees the existence of some $\,b\in R\,$ such that
$\,a+e\,b=a+t\,(y\,b)\in {\rm U}(R)$. This checks that
$\,{\rm sr}\,(a)=1$.\qed

\bs
Coming back to Theorem 2.2 itself, we collect in the next result
a number of other basic facts on the notion of ``element-wise
right stable range one''.

\bs\nt
{\bf Theorem 2.4.} (A) {\it For any ring element $\,a\in R$,
$\,{\rm sr}\,(a)=1\,$ iff, for any $\,t\in R$, $\,aR+tR=R\,$
implies that $\,{\rm sr}\,(a+ty)=1\,$ for some $\,y\in R$.}\\
(B) {\it For any $\,a\in R\,$ and $\,u,v\in {\rm U}(R)$,
$\,{\rm sr}\,(a)=1\,\Rightarrow\,{\rm sr}\,(uav)=1$.}\\
(C) {\it If $\,a\,$ has a one-sided inverse in $\,R$, then
$\,{\rm sr}\,(a)=1\,$ iff $\,a\in {\rm U}(R)$.}\\
(D) {\it For any ideal $\,J\subseteq R$, let $\,\overline{R}
=R/J$. For any $\,a\in R$, $\,{\rm sr}\,(a)=1\,\Rightarrow\,
{\rm sr}\,\bigl(\overline{a}\bigr)=1$. The converse holds if the
quotient map $\,f:R\rightarrow R/J\,$} ``{\it reflects units\/}'',
{\it in the sense that $\,f(a)\in {\rm U}(R/J) \Rightarrow a\in
{\rm U}(R)$. In particular, if $\,J\subseteq {\rm rad}\,(R)\,$}
({\it the Jacobson radical of $\,R\,$}), {\it then $\,{\rm sr}_R(a)
=1\,$ iff $\,{\rm sr}_{R/J}\bigl(\overline{a}\bigr)=1$.}

\bs\nt
{\bf Proof.} (A) The ``only if" part is trivial as we can
choose $\,y=0$. To prove the ``if\/" part, consider any
equation $\,aR+tR=R$.  By assumption, there exists $\,y\in
R\,$ such that $\,a'=a+ty\,$ has $\,{\rm sr}\,(a')=1$. Then
$\,a'R+tR\supseteq aR+tR=R\,$ implies that there exists $\,y_0
\in R\,$ such that $\,a'+ty_0\in {\rm U}(R)$. As $\,a'+ty_0
=a+t\,(y+y_0)$, this checks that $\,{\rm sr}\,(a)=1$.


\mk\nt
(B) After replacing $\,uav\,$ by its conjugate $\,u^{-1}(uav)\,u
=avu$, we need only show that $\,{\rm sr}\,(av)=1$. If $\,avR+tR=R$,
we have $\,aR+tR=R$, so $\,a+tb\in {\rm U}(R)\,$ for some $\,b\in R$.
This implies that $\,av+t\,(bv)\in {\rm U}(R)$, which checks the claim
that $\,{\rm sr}\,(av)=1$.

\mk\nt
(C) The result in (C) here should be regarded as an ``element-wise
version'' of the well known fact that rings of stable range one (see
(2.5A) below) are always Dedekind-finite. In view of what we said
at the beginning of Example 2.1, we need only prove the ``only if''
part in (C), so assume that $\,{\rm sr}\,(a)=1$. If $\,ax=1\,$ for
some $\,x\in R$, applying (2.2)(5) for $\,c=0\,$ shows directly
that $\,a\in {\rm U}(R)$.  If $\,ya=1\,$ for some $\,y\in R\,$
instead, applying (2.2)(4) yields some $\,b\in R\,$ such that
$\,u\!:=a+(1-ay)\,b\in {\rm U}(R)$. Left multiplying this equation
by $\,y\,$ shows that $\,yu=ya+(y-yay)\,b=1$.  Thus, $\,y\in
{\rm U}(R)$, and hence $\,a\in {\rm U}(R)$.


\mk\nt
(D) The two main conclusions in (D) follow from a routine application
of the criterion for stable range one elements given in Theorem 2.2(3).
Given these two conclusions, the last part of (D) follows from the
well known fact that the quotient map $\,R \rightarrow R/J\,$
reflects units if $\,J\subseteq {\rm rad}\,(R)$.\qed




\bs\nt
{\bf Example 2.5A.} Recall that a ring $\,R\,$ has stable range one
iff $\,{\rm sr}_R(a)=1\,$ for every $\,a\in R$.  Well known classes
of such rings are: (1) semilocal rings (including local rings and
one-sided artinian rings according to [La$_{1}$:~(20.3)] and
[La$_3$:~(2.10)]); (2) unit-regular rings (including all abelian
regular rings) according to [Go:~(4.12)]; (3) strongly $\,\pi$-regular
rings according to [Ar:~Theorem 4] (see also [Khu]); (4) left-and-right
self-injective rings according to [La$_3$: Cor.~7.17]; and (5) the ring
of all algebraic integers according to [Va$_2$:~(1.2)].  Many other
interesting examples of commutative and noncommutative rings of stable
range one were given in [Va$_1$, Va$_2$], [EO], [GM], [CY$_1$, CY$_2$],
[Yu$_1$, Yu$_2$], [Ch$_3$, Ch$_4$], [La$_3$], etc.

\bs\nt
{\bf Example 2.5B.} Using (2.2) and (2.4), it is easy to determine
which integers $\,n\,$ have stable range one in the ring
$\,{\mathbb Z}$. The answer is, $\,{\rm sr}\,(a)=1\,$ iff $\,a\in
\{0,\,{\pm 1}\}$. The ``if'' part is clear (say from Example 2.1).
For the ``only if'' part, start with $\,{\rm sr}\,(n)=1\,$ for some
$\,n\in {\mathbb Z}$.  By Theorem 2.4(B), we may replace $\,n\,$ by
$\,-n\,$ if necessary to assume that $\,n\ge 0$. If $\,n\geq 2$,
taking $\,x=-n\,$ in Theorem 2.2(3) shows that there exists $\,b\in
{\mathbb Z}\,$ such that $\,n+(1+n^2)\,b\in \{{\pm 1}\}$.  This is
clearly impossible, so we must have $\,n\in \{0,1\}$.  Note that
this example shows that the set of stable range one elements in
a ring is in general not closed under the map $\,a\mapsto 1-a\,$.
It shows also that a ``clean element'' (to be defined in Example 2.5F
below) such as $\,2\in R={\mathbb Z}\,$ need not have stable range
one in $\,R$.

\bs\nt
{\bf Example 2.5C.} Since units and idempotents in any ring
always have stable range one, it will follow from the forthcoming
Product Theorem 2.8 that {\it unit-regular elements also have stable
range one.}  Somewhat surprisingly, however, in any matrix ring
$\,R={\mathbb M}_n(S)\,$ with $\,n\geq 2$, all matrices $\,aE_{ij}\,$
with $\,a\in S\,$ (i.e.~matrices with at most one nonzero entry)
happen to have also stable range one. This interesting fact will
be proved later in Theorem 5.8.

\bs\nt
{\bf Example 2.5D.} Some fairly large rings may have ``very few''
elements of stable range one. For instance, let $\,R\,$ be a polynomial
ring $\,S\,[x]\,$ over a commutative domain $\,S$. Then the set of
stable range one elements in $\,R\,$ is just $\,\{0\}\cup {\rm U}(S)\,$
(although $\,S\,$ itself may have many other elements of stable range
one). Indeed, if a degree $\,n\,$ polynomial $\,a(x)\in R\setminus
\{0\}\,$ has stable range one, then by Theorem 2.2, there exists
$\,b(x)\in R\,$ such that
$$
f(x)\!:=a(x)+(1-a(x)\,x)\,b(x)\in {\rm U}(R)={\rm U}(S).
$$
If $\,b(x)\neq 0$, we have $\,{\rm deg}\,f(x)\geq n+1$, which
is impossible. Thus, we must have $\,b(x)=0$, and so
$\,f(x)=a(x)\in {\rm U}(S)$.

\bs\nt
{\bf Example 2.5E.} For an example of a ring with a relatively large
set of stable range one elements, let $\,R=C(X)\,$ be the ring of
continuous real-valued functions on a topological space $\,X$. For
any $\,a\in R\,$ such that $\,a\,(X)\geq 0$, Azarpanah, Farokhpay
and Ghashghaei have observed in [AFG:~(3.2)] that $\,{\rm sr}\,(a)=1$.
Indeed, if $\,b\in R\,$ is such that $\,aR+bR=R$, then $\,a,\,b\in
C(X)\,$ cannot have a common zero. This implies that the continuous
function $\,a+b^2\,$ takes only {\it positive\/} values on $\,X$,
and so $\,a+b^2\in {\rm U}(R)$.  In the case where $\,X\,$ is a
completely regular space, a precise description of the set of stable
range one elements in $\,C(X)\,$ has been given by Ghashghaei [Gh].

\bs\nt
{\bf Example 2.5F.} There are certainly some connections between
stable range one elements and Nicholson's notion of {\bf clean
elements\/} in rings.  In Nicholson's paper [Ni$_1$], a ring
element $\,a\in R\,$ is called {\it clean\/} if $\,a\,$ is the sum
of an idempotent and a unit in $\,R$. In view of Example 2.1, such
an element $\,a\,$ is the sum of two elements of stable range one
in $\,R$, although $\,a\,$ itself may not have stable range one.
Going in the other direction, consider for instance any two
elements $\,a,b \in R\,$ such that $\,aR+bR=R$.  If either
$\,{\rm sr}_R(a)=1\,$ or $\,{\rm sr}_R(b)=1$, it turns out that
$\,\footnotesize{\begin{pmatrix}a&b\\0&0\end{pmatrix}}$ is a
clean matrix in $\,{\mathbb M}_2(R)$.  For a proof of this
interesting fact, based on the ideas developed in [KL$_1$:~\S3],
see ``Example (11)'' listed under (17.2) in [La$_4$].

\bs
In the next three theorems, we shall give some more applications of
Theorem 2.4.  The first one is the following basic result relating
the Jacobson radical $\,{\rm rad}\,(R)\,$ of a ring $\,R\,$ to the
element-wise notion of (right) stable range one.

\bs\nt
{\bf Theorem 2.6.}
(1) {\it For any $\,a\in R\,$ and $\,b\in {\rm rad}\,(R)$,
we have $\,{\rm sr}\,(a)=1\,$ iff $\,{\rm sr}\,(a+b)=1$.}
({\it In particular, $\,{\rm sr}\,(b)=1\,$ always holds.})\\
(2) {\it For any $\,b\in R$, we have $\,b\in {\rm rad}\,(R)\,$ iff
$\,{\rm sr}\,(b)=1\,$ and $\,b+{\rm U}(R)\subseteq {\rm U}(R)$.}

\bs\nt
{\bf Proof.} (1) This part follows quickly from Theorem 2.4(D)
since the natural surjection $\,R\rightarrow R/{\rm rad}\,(R)\,$
reflects units, and both conditions in (1) are equivalent to
$\,{\rm sr}\,\bigl(\overline{a}\bigr)=1\,$ in the factor ring
$\,\overline{R}=R/{\rm rad}\,(R)$.


\mk\nt
(2) The ``only if'' part follows from (1) by setting $\,a=0$,
and from a familiar property of the Jacobson radical.
The ``if'' part is basically folklore after the appearance of
[La$_2$:~Exercise 20.10B].  For more details, assume that
$\,{\rm sr}\,(b)=1\,$ and $\,b+{\rm U}(R)\subseteq {\rm U}(R)$.
To get the desired conclusion that $\,b\in {\rm rad}\,(R)$, it
suffices to check that, for any $\,x\in R$, $\,1-bx\,$ has a right
inverse in $\,R$. (Here, we are using a classical characterization
result for elements in $\,{\rm rad}\,(R)$, as given for example in
[La$_1$:~Lemma 4.1].) By Theorem 2.2, there exists $\,y\in R\,$
such that $\,b+(1-bx)\,y\in {\rm U}(R)$.  Thus, $\,(1-bx)\,(-y)
\in b+{\rm U}(R)\subseteq {\rm U}(R)$, which shows that $\,1-bx\,$
has a right inverse, as desired.\qed


\bs
Our next application of Theorem 2.4 is largely inspired by the
work of Godefroid [Gd] and Goodearl-Menal [GM].


\bs\nt
{\bf Theorem 2.7.} {\it Suppose $\,a\in R\,$ has the property that,
for any $\,x\in R$, there exists $\,u\in {\rm U}(R)\,$ such that
$\,x-u^{-1}\in {\rm U}(R)\,$ and $\,{\rm sr}\,(a-u)=1$.  Then
$\,{\rm sr}\,(a)=1$.}

\bs\nt
{\bf Proof.} It suffices to check the validity of the condition (A)
in Theorem 2.4, so consider any equation $\,ax+ty=1$. By assumption,
there exists $\,u\in {\rm U}(R)\,$ such that $\,v\!:=x-u^{-1}
\in {\rm U}(R)\,$ and $\,{\rm sr}\,(a-u)=1$. We have then
$$
av+ty=a\,(x-u^{-1})+ty=1-au^{-1}=(a-u)\,(-u^{-1}),
$$
so $\,a+t\,(yv^{-1})=(a-u)\,(-u^{-1}v^{-1})$. According to
Theorem 2.4(B), this element has stable range one, so we
have checked that the condition in Theorem 2.4(A) holds.\qed

\bs
We note in passing that Theorem 2.7 applies somewhat less
widely than Theorem 2.2 and Theorem 2.4.  For instance,
it cannot be applied to a ring $\,R\,$ unless every
element in $\,R\,$ is a sum of two units.

\bs
Applying Theorem 2.4 also retrieves easily the following result
that is due to H.~Chen and W.\,K.~Nicholson; see [CN:~Lemma 17].
This important result will be used over and over again in the
forthcoming sections.


\bs\nt
{\bf Product Theorem 2.8.}  {\it If $\,{\rm sr}\,(a_i)=1\,$ for
all $\,i\in [1,n]$, then $\,{\rm sr}\,(a_1\cdots a_n)=1$.\/}
({\it Equivalently, the set} $\,\{a\in R:\,{\rm sr}\,(a)=1\}\,$
{\it is a monoid under multiplication\/} ({\it with identity
$\,1\,$}) {\it containing $\,{\rm U}(R)\,$ as a subgroup.})

\bs\nt
{\bf Proof.} It suffices to check that $\,{\rm sr}\,(a)
={\rm sr}\,(a')=1\Rightarrow {\rm sr}\,(aa')=1$.  Following [CN], 
we start with an equation 
$\,aa'R+tR=R$. Then $\,aR+tR=R\,$ too, so $\,{\rm sr}\,(a)
=1\,$ implies that $\,u\!:=a+t\,b\in {\rm U}(R)\,$ for
some $\,b\in R$.  Thus, $\,aa'+t\,(ba')=ua'$, with
$\,{\rm sr}\,(ua')=1\,$ by Theorem 2.4(B).  According to
Theorem 2.4(A), this checks that $\,{\rm sr}\,(aa')=1$,
as desired. \qed

\bs\nt
{\bf Example 2.9.} In contrast with the theorem above, the well known
example of $\,R={\mathbb Z}\,$ shows that the set of stable range one
elements in a ring $\,R\,$ may not be closed under addition. However,
one may still ask if $\,{\rm sr}\,(a)={\rm sr}\,(b)=1\,$ might imply
$\,{\rm sr}\,(a+b)=1\,$ if we assume further that $\,ab=ba=0$.
Unfortunately, this additional condition won't do it either. Indeed,
in the ring $\,R={\mathbb M}_2\bigl({\mathbb Z}\bigr)$, both of the
matrices $\,A={\rm diag}\,(2,0)\,$ and $\,B={\rm diag}\,(0,2)\,$ have
stable range one by our later result (5.8) in \S5, and they satisfy
the ``orthogonality condition'' $\,AB=BA=0$.  However, $\,A+B=2I_2\,$
does not have stable range one by our later result (7.5) in \S7.


\bs
To conclude this section, we rcord the following easy
consequence of Theorem 2.8.

\bs\nt
{\bf Corollary 2.10.} {\it If $\,xay=1\in R\,$ where
$\,{\rm sr}\,(a)=1\,$ and either $\,{\rm sr}\,(x)=1\,$
or $\,{\rm sr}\,(y)=1$, then $\,a\in {\rm U}(R)$.}

\bs\nt
{\bf Proof.} We may assume that $\,{\rm sr}\,(x)=1\,$ (since the
case where $\,{\rm sr}\,(y)=1\,$ is similar). Here, Theorem 2.8
gives $\,{\rm sr}\,(xa)=1$. Since $\,xa\,$ has a right inverse $\,y$,
Theorem 2.4(C) shows that $\,z\,(xa)=1\,$ for some $\,z\in R$.
Thus, $\,a\,$ has a left inverse $\,zx$, so another application
of Theorem 2.4(C) shows that $\,a\in {\rm U}(R)$. \qed

\mk

\bs\nt
{\bf \S3. \ Super Jacobson's Lemma and Its Applications}

\bs
One main question in the study of element-wise stable range one
in arbitrary (noncommutative) rings is whether such a notion is
always left-right symmetric for all elements. In the literature
to date, this question has been answered positively only for a
few specific classes of rings; e.g.~regular rings (as in \S4),
duo rings\footnote{A ring $\,R\,$ is said to be a {\bf duo ring\/}
if every 1-sided ideal in $\,R\,$ is a 2-sided ideal.} (as in [AO$_2$]),
and the ring $\,{\mathbb M}_2\bigl({\mathbb Z}\bigr)\,$ (as
in [CP$_1$]). The first result in this section is devoted to
resolving the question above positively {\it for all rings.}

\bs\nt
{\bf Symmetry Theorem 3.1.} {\it For any ring element $\,a\in R$,
the property $\,{\rm sr}\,(a)=1\,$ is {\bf equivalent\/} to the
property $\,{\rm sr}\,'(a)=1$. Here, as in {\rm Definition 1.1,}
$\,{\rm sr}\,'(a)=1\,$ refers to the property that $\,a\in R\,$ has
{\bf left\/} stable range one, while $\,{\rm sr}\,(a)=1\,$ means,
as always, that $\,a\in R\,$ has {\bf right\/} stable range one.}

\bs
As it turned out, the main trick for the proof of the above symmetry
result is the following interesting fact in noncommutative ring
theory, which we shall fondly call by a special name in honor of
the late Yale Professor Nathan Jacobson.  Quite remarkably, this
result applies uniformly to all three of the sets $\,{\rm U}(R)\,$
(the group of units), $\,{\rm reg}\,(R)\,$ (the set of regular
elements), and $\,{\rm ureg}\,(R)\,$ (the set of unit-regular
elements), although not to $\,{\rm sreg}\,(R)\,$ (the set of
strongly regular elements).

\bs\nt
{\bf Super Jacobson's Lemma 3.2.} {\it For any three                           elements $\,a,b,x\,$ in any ring $\,R$,\break
the following three conclusions hold\,}:

\mk\nt
(1) {\it $\,a+b-axb\in {\rm U}(R)\,$ iff $\,a+b-bxa\in {\rm U}(R)$.}\\
(2) {\it $\,a+b-axb\in {\rm reg}(R)\,$ iff $\,a+b-bxa\in {\rm reg}(R)$.}\\
(3) {\it $\,a+b-axb\in {\rm ureg}(R)\,$ iff $\,a+b-bxa\in {\rm ureg}(R)$.}

\mk\nt
{\it However, in general, $\,a+b-axb\in {\rm sreg}(R)\,$ is not
equivalent to $\,a+b-bxa\in {\rm sreg}(R)$.}


\bs
After the first draft of this paper was circulated, Professor P.~Ara
has kindly informed the authors that Case (1) above was first proved
by Menal and Moncasi in [MM: p.\,296].  As soon as Case (1) of this
lemma is granted, a quick recall of the criterion (2) in Theorem 2.2
for $\,{\rm sr}\,(a)=1\,$ and the corresponding criterion for
$\,{\rm sr}\,'(a)=1\,$ will show that these two properties are
{\it equivalent\/}, for any element $\,a\in R\,$ in any ring $\,R$,
thus proving the Symmetry Theorem 3.1.  In the language introduced
by C\v{a}lug\v{a}reanu and Pop in [CP$_1$], if $\,{\rm sr}\,(a)=1\,$
for a given element $\,a\in R$, then for any $\,x\in R$, any
``right unitizer'' $\,b\,$ for $\,x\,$ will automatically be
a ``left unitizer'' for $\,x$.  This may sound almost too good
to be true, but the following proof will show that it is\,!


\bs\nt
{\bf Proof of Lemma 3.2.} To prove this ``Super Lemma'', we will
first focus on Case (1); namely, the case of the unit group
$\,{\rm U}(R)$. Of course, it will be sufficient to prove the
``only if'' part of the lemma.  Before we proceed formally with
the proof, perhaps a short explanation on the nomenclature used
in this lemma would be helpful.  Consider the special case where
$\,x=1$, and assume that $\,a+b-ab\in {\rm U}(R)$. Then
$$
1-(1-a)\,(1-b)=1-(1-a-b+ab)=a+b-ab\in {\rm U}(R).
$$
By the classical Jacobson's Lemma (as recalled in the Abstract),
this implies that
$$
1-(1-b)\,(1-a)=1-(1-a-b+b\,a)=a+b-b\,a\in {\rm U}(R),
$$
which is what we want. Thus, to prove Lemma 3.2 in the Case (1), our
job is to come up with a generalization of the above implication when
a {\it general\/} element $\,x\in R\,$ is given such that $\,a+b-axb
\in {\rm U}(R)$. To accomplish this goal, we will exploit the use of
$\,2\times 2\,$ matrices over our ring $\,R$. Although we can no longer
use the tool of ``determinants'' for such matrices, it turns out that a
suitable appeal to the basic idea behind the proof of ``Banachiewicz's
Inversion Theorem'' in [Ban] comes to our rescue.  In its ``modern''
{\it purely ring-theoretic\/} incarnation as given in [La$_4$:~(10.25)],
this well known result states that, for any $\,u\in {\rm U}(R)$,
a matrix $\,A=\footnotesize{\begin{pmatrix} u&q\\p&r\end{pmatrix}}
\in {\mathbb M}_2(R)\,$ is invertible iff $\,r-pu^{-1}q\in
{\rm U}(R)$.\footnote{The detailed expression for $\,A^{-1}\,$ in
this case is given explicitly in [La$_4$:~(10.29)].  However, this
expression is not needed for any of the results in this section.}
In our special setting, we may start with the following matrix-product
equation which has appeared before, for instance, in [LN:~\S4]:
$$
\begin{pmatrix}1&-b\\1-ax&a\end{pmatrix}
  \begin{pmatrix}1&0\\x&1\end{pmatrix}
= \begin{pmatrix}1-b\,x&-b\\1&a\end{pmatrix}
      \in {\mathbb M}_2(R). \leqno (3.3)
$$
If we assume that $\,a+b-axb\in {\rm U}(R)$, then Banachiewicz's
Inversion Theorem implies that the first matrix in (3.3) above
is invertible, so the last matrix there is also invertible.
After a row permutation, it follows that
$\footnotesize{\begin{pmatrix}1&a\\1-b\,x&-b\end{pmatrix}}$ is
invertible.  Applying Banachiewicz's Inversion Theorem  once
more shows that $\,a+b-bxa\in {\rm U}(R)$, as desired.

\mk
Having done the work above, we need to handle two more cases; namely,
the cases where $\,{\rm U}(R)\,$ is replaced by $\,{\rm reg}\,(R)\,$
or by $\,{\rm ureg}\,(R)$. Fortunately, not much further work is
needed to get the desired conclusions in these two cases. Indeed, in
the relatively new treatment of Banachiewicz's inversion result in
[La$_4$:~Theorem 10.25], it is already proved that this result also
works in the cases where $\,A=\footnotesize{\begin{pmatrix}u&q\\p&r
\end{pmatrix}}$ is regular or unit-regular in $\,{\mathbb M}_2(R)$,
with the respective characterizations that $\,r-pu^{-1}q\in
{\rm reg}\,(R)\,$ or $\,r-pu^{-1}q\in {\rm ureg}\,(R)$. This proves
the truth of the desired conclusions in both Case (2) and Case (3).

\mk
To show that the above work cannot be extended to the case of
$\,{\rm sreg}\,(R)$, we need only consider the following example.
In the ring $\,R={\mathbb M}_2\bigl({\mathbb Z}\bigr)$, let
$\,a=\footnotesize{\begin{pmatrix}1&1\\0&0\end{pmatrix}}$,
$\,b=\footnotesize{\begin{pmatrix}1&0\\0&0\end{pmatrix}}$, and
$\,x=\footnotesize{\begin{pmatrix}0&1\\1&0\end{pmatrix}}$. Then
$\,a+b-axb=\footnotesize{\begin{pmatrix}1&1\\0&0\end{pmatrix}}
\in {\rm idem}\,(R)\subseteq {\rm sreg}\,(R)$. However,
$\,c\!:=a+b-bxa=\footnotesize{\begin{pmatrix}2&1\\0&0
\end{pmatrix}}$ is clearly {\it not\/} in $\,R\cdot c^2\,$ since
$\,c^2=\footnotesize{\begin{pmatrix}4&2\\0&0\end{pmatrix}}$.
This shows that $\,c\notin {\rm sreg}\,(R)$. (On the other hand,
we do have $\,c\in {\rm ureg}\,(R)$, so there is no contradiction
to the conclusions obtained in Case (2) or Case (3).)  \qed


\bs\nt
{\bf Example 3.4.} Talking about ``ternary generalizations'' of
Jacobson's Lemma, we should point out that other ``more obvious''
(or ``better looking'') forms of such generalizations may not hold
true at all. For instance, given arbitrary $\,a,x,b\in R$, one may
wonder if $\,1-axb\in {\rm U}(R)\,$ might be equivalent to $\,1-bxa
\in {\rm U}(R)$.  The answer to this question is, unfortunately,
``no'' (in general), as the following easy example shows.
In any nonzero matrix ring $\,R={\mathbb M}_2(S)$, let
$\,a=\footnotesize{\begin{pmatrix}1&0\\0&0\end{pmatrix}}$,
$\,x=\footnotesize{\begin{pmatrix}0&1\\1&0\end{pmatrix}}$,
and $\,b=\footnotesize{\begin{pmatrix}1&1\\0&0\end{pmatrix}}$.
Then $\,1-axb=1\in {\rm U}(R)$, but
$\,1-bxa=\footnotesize{\begin{pmatrix}0&0\\0&1\end{pmatrix}}
\notin {\rm U}(R)\,$!  On the other hand, for the same
choices of $\,a,b,x\in R$, it is clear that neither
$\,a+b-axb=\footnotesize{\begin{pmatrix}2&1\\0&0\end{pmatrix}}$
nor $\,a+b-bxa=\footnotesize{\begin{pmatrix}1&1\\0&0\end{pmatrix}}$
is in $\,{\rm U}(R)$, thus reaffirming the truth of Super Jacobson's
Lemma in this particular case.


\bs\nt
{\bf Historical Remark 3.5.} As we have seen from the calculations
given at the beginning of the proof of (3.2), the classical form
of Jacobson's Lemma corresponds to the special case $\,x=1\,$ of
Super Jacobson's Lemma when we make the substitutions $\,a\mapsto
1-a\,$ and $\,b\mapsto 1-b$. The new formulation of Super Jacobson's
Lemma 3.2 also suggests that we can generalize another old idea of
Jacobson, as follows. In Jacobson's book [Ja$_2$:~p.8], an associative
``circle operation'' for ring elements $\,a,b\in R\,$ was introduced
through the definition $\,a\circ b :=a+b-ab$. If we pick a
{\bf fixed\/} element $\,x\in R$, we can define a new circle
operation by the rule $\,a\circ b\!:=a+b-axb$.  Quite pleasantly,
this still turns out to be an associative operation, since a quick
calculation shows that $\,(a\circ b) \circ c\,$ and $\,a\circ
(b\circ c)\,$ are both equal to $\,a+b+c-axb-bxc-axc+axbxc$.
Furthermore, $0\,$ is the identity element for this new circle
operation; that is, $\,a\circ 0=a=0\circ a\,$ for every
$\,a\in R$.\footnote{Note that, in general, this new operation is
not commutative. Indeed, we can check easily that it is commutative
iff $\,x\,$ is a central element of $\,R\,$ that annihilates all
``additive commutators'' $\,ab-ba\,$ in $\,R$.} Expressed in terms
of the new circle operation (defined via a given element $\,x\,$),
Super Jacobson's Lemma would simply say that $\,a\circ b\in
{\rm U}(R)\,$ iff $\,b\circ a\in {\rm U}(R)$; $\,a\circ b\in
{\rm reg}\,(R)\,$ iff $\,b\circ a\in {\rm reg}\,(R)$; and
$\,a\circ b\in {\rm ureg}(R)\,$ iff $\,b\circ a\in {\rm ureg}(R)$.

\bs
To further illustrate the power of Super Jacobson's Lemma 3.2,
we shall next state and prove below four of its interesting
``binary specializations'' for any ring $\,R$.

\bs\nt
{\bf Proposition 3.6.} {\it For any $\,a,x\in R$, the following
conditions are equivalent\/}:

\sk\nt
(1) $\,1-ax+axa\in {\rm U}(R)$.\;\;\; (2) $\,1-xa+axa\in {\rm U}(R)$.\\
(3) $\,1-ax+a^2 x\in {\rm U}(R)$.\;\;\;\, (4) $\,1-xa+x a^2\in {\rm U}(R)$.

\mk\nt
{\it If these conditions hold, then $\,\{ax-xa,\;axa-xa^2,\;axa-a^2x\}
\subseteq {\rm U}(R)+{\rm U}(R)$. Also, the same result holds
if $\,{\rm U}(R)\,$ is replaced throughout by $\,{\rm reg}\,(R)$,
or by $\,{\rm ureg}\,(R)$.}

\bs\nt
{\bf Proof.} It suffices to work with the case of $\,{\rm U}(R)$, since
the other two cases (for $\,{\rm reg}\,(R)\,$ and $\,{\rm ureg}\,(R)$)
are completely similar. To begin with, $(1)\Leftrightarrow (2)$
follows from Super Jacobson's Lemma by taking $\,b=1-a$.  Next,
$(1)\Leftrightarrow (4)$ follows from the {\it usual\/} Jacobson's
Lemma (in all three cases) by writing $\,1-ax+axa=1-a\,(x-xa)\,$ and
$\,1-xa+xa^2=1-(x-xa)\,a$.  With this little trick at our disposal,
$(2)\Leftrightarrow (3)$ follows similarly. Finally, the conclusion
about the three elements $\,ax-xa,\;axa-xa^2\,$ and $\,axa-a^2x\,$
follows by subtracting (1) from (2); (4) from (2); and (3) from (1),
respectively.\qed

\bs
Comparing Proposition 3.6 with what is known in the literature,
we can say the following.  In the special case where we assume
additionally that $\,a=axa$, (1) and (2) were known to be both
equivalent to $\,a\in {\rm sreg}\,(R)\,$ (i.e.~$a\,$ is
{\it strongly regular\/} in $\,R\,$) by a result of Puystjens and
Hartwig in [PuH], while (3) and (4) can be easily shown to be
equivalent to (1) and (2).  For a detailed exposition on these
four characterizations for the elements in $\,{\rm sreg}\,(R)$,
see [La$_4$:~Theorem 3.24]. The remarkable thing about Proposition
3.6 is that it holds {\it without any initial binding conditions\/}
between the two elements $\,a,\,x\in R$, while it generalizes the
work of Puystjens and Hartwig on strongly regular elements in
a rather unexpected way.

\bs
Upon reviewing the meaning of Super Jacobson's Lemma in the special case
where $\,R={\mathbb M}_n(S)\,$ for some ring $\,S$, we see that, for
any three matrices $\,A,B,X\in R$, we have $\,A+B-AXB\in {\rm GL}_n(S)\,$
iff $\,A+B-BXA\in {\rm GL}_n(S)$. In the special case where $\,S\,$ is
a {\it commutative\/} ring, it is well known that a matrix $\,M\in R\,$
is in $\,{\rm GL}_n(S)\,$ iff $\,{\rm det}\,(M)\in {\rm U}(S)$. In this
special case, Super Jacobson's Lemma for the matrix ring $\,R\,$ above
would amount to the statement that $\,{\rm det}\,(A+B-AXB)\in {\rm U}(S)\,$
iff $\,{\rm det}\,(A+B-BXA)\in {\rm U}(S)$. This was the first time it
occurred to us that, for any $\,A,B,X \in R$, the two determinants
above must be rather ``intimately related''.  Indeed, after some
additional work, the following result about determinants was proved
in [KL$_3$: Theorem 2.1].

\bs\nt
{\bf Theorem 3.7.} {\it For any $\,A,B,X\in R={\mathbb M}_n(S)\,$
over a commutative ring $\,S$, we have}
$$
{\rm det}\,(A+B-AXB)={\rm det}\,(A+B-BXA)\in S. \leqno (3.8)
$$
{\it In particular, $\,A+B-AXB \in {\rm GL}_n(S)\,$ if and only
if $\,A+B-BXA \in {\rm GL}_n(S)$.  However, we may not have
$\,{\rm tr}\,(A+B-AXB)={\rm tr}\,(A+B-BXA)$.}

\bs
While we will not repeat the proof of the equation (3.8) in [KL$_3$]
here, we should point out that, in studying the possible applications
of the Symmetry Theorem 3.1, it is perhaps to be expected that the
notion of a \starring\, (or a ``ring with involution'') would be of
interest here.  Indeed, using the Symmetry Theorem 3.1, we can prove
easily the following result on ``element-wise stable range one'' for
elements in any \starring.

\bs\nt
{\bf Theorem 3.9.} {\it If \Rstar\, is a \starring\,, then
for any $\,a\in R$, $\,{\rm sr}\,(a)=1\,$ iff $\,{\rm sr}\,(a^*)=1$.}

\bs\nt
{\bf Proof.} Thanks to the Symmetry Theorem 3.1, it suffices to
prove that, if $\,{\rm sr}\,(a)=1$, then $\,{\rm sr}\,'(a^*)=1$.
To do this, we take any element of $\,R\,$ (which we may
conveniently write as $\,x^*\,$ for some $\,x\in R$\/), and try
to show that $\,a^*+b^*\,(1-x^* a^*)\in {\rm U}(R)\,$ for some
$\,b\in R$. Since $\,{\rm sr}\,(a)=1$, there exists $\,b\in R\,$
such that $\,a+(1-ax)\,b\in {\rm U}(R)$. Applying the involution
$\,^*\,$ to this membership relation (and noting that $^*$
preserves the unit group $\,{\rm U}(R)$), we see that
$\,a^*+b^*\,(1-x^* a^*)\in {\rm U}(R)$, as desired. \qed

\bs
Over a {\it commutative\/} ground ring, a quick application
of Theorem 3.9 yields the following nice result on the
transposition of square matrices over a commutative ring.

\bs\nt
{\bf Corollary 3.10.} {\it If $\,R={\mathbb M}_n(S)\,$ where
$\,S\,$ is any commutative ring, then for any matrix $\,A\in R$,
$\,{\rm sr}\,(A)=1\,$ iff $\,{\rm sr}\,(A^T)=1$, where $\,A^T\,$
denotes the transpose of $\,A$.}

\bs\nt
{\bf Proof.} We simply note that, in this case, the transpose
operation $\,A\mapsto A^T\,$ is a ring involution on
$\,R={\mathbb M}_n(S)$, and apply Theorem 3.9.\qed

\bs
While the proof given above definitely made use of the commutativity
property of the base ring $\,S$, it is nevertheless legitimate to ask
whether the conclusion of Corollary 3.10 might remain true if $\,S\,$
is just any ring. We will return to answer this question negatively
later; see Example 6.12.

\bs\nt
{\bf Remark 3.11.} To conclude this section, we would like to make
the following further observations about the technique of applying
Super Jacobson's Lemma to the study of the notion of stable range one.
In the standard literature, it is known that this notion is also worth
studying in some of its variations, such as {\it unit stable range one\/},
{\it idempotent stable range one\/}, {\it regular stable range one\/},
and {\it square stable range one\/}; see, for instance, Goodearl-Menal
[GM], Chen [Ch$_1$], Ashrafi-Nasibi [AN], and Wang et al.~[WC], [KLW].
These notions arise naturally by specializing the ``element-wise stable
range one'' notion in Definition 1.1 by asking that the existential
element $\,b\in R\,$ be, respectively, a unit, an idempotent, a regular
element, or a perfect square in $\,R$.  While these somewhat more
refined notions have again a ``left version'' and a ``right version''
(depending on the use of left ideals or right ideals in $\,R\,$), the
very same argument used for proving the Symmetry Theorem 3.1 can be
exploited once more to show that they are {\it all\/} left-right
symmetric notions as well.

\bs\nt
{\bf Example 3.12.} Of course, it is entirely to be expected that the
more refined notions of stable range one mentioned in Remark 3.11
above are in general more restrictive than the original notion of
stable range one introduced in Definition 1.1. For some quick examples
to show some cases of this, consider the matrix $\,A={\rm diag}\,(7,0)
\in R={\mathbb M_2}\bigl({\mathbb Z}\bigr)$, with $\,{\rm sr}\,(A)
=1\,$ to be shown in Theorem 5.8. We claim that, in $\,R$,
$\,A\,$ has neither unit stable range one nor idempotent stable
range one.  To see this, let $\,B={\rm diag}\,(2,1)\,$ and
$\,C={\rm diag}\,(-3,1)$. Since $\,A+B\,C=I_2$, we have $\,AR+BR=R$.
If $\,A\,$ has unit stable range one, there would exist a matrix
$\,U=\footnotesize{\begin{pmatrix}a&b\\c&d\end{pmatrix}}\in
{\rm U}(R)\,$ such that $\,A+B\,U\in {\rm U}(R)$. A quick computation
shows that
$$
{\rm det}\,(A+B\,U)=7d+2\,(ad-bc)=7d\,{\pm}\,2\notin \{{\pm 1}\},
$$
which is a contradiction.  Somewhat similarly, if $\,A\,$
has idempotent stable range one, there would exist
$\,P=\footnotesize{\begin{pmatrix}p&q\\r&s\end{pmatrix}}\in
{\rm idem}\,(R)\,$ such that $\,A+BP\in {\rm U}(R)$. Clearly,
$\,P\notin \{0,\,I_2\}$, so we must have $\,p+s=1\,$ and $\,ps-qr=0$.
A quick computation shows that $\,{\rm det}\,(A+B\,P)=7\,(1-p)
\notin \{\pm 1\}$, again a contradiction.

\mk

\bs\nt
{\bf \S4. \ Stable Range One Elements Via Unit-Regularity}

\bs
To further develop the notion of ring elements of stable range one,
we go back to the basic theme of looking at regular and unit-regular
elements in any ring. The first theorem in this section is largely
motivated by a classical result of Fuchs, Kaplansky and Henriksen
on regular rings that is reported in Goodearl's book [Go:~(4.12)].
According to this result, a (von Neumann) regular ring has
stable range one iff it is unit-regular.  To make this into an
``element-wise result'', we have the following theorem characterizing
regular elements with stable range one {\it in any ring.}  Here,
$\,(1)\Leftrightarrow (2)$ was noted in [KL$_2$: (3.5)] (see also
[HN: (3.2)]), while $\,(2)\Leftrightarrow (4)$ was an observation
of H.~Chen reported in [HN: (5.13)].



\bs\nt
{\bf Theorem 4.1.} {\it For any element $\,a\in {\rm reg}\,(R)$,
the following conditions are equivalent\,}:

\mk\nt
(1) $\,{\rm sr}\,(a)=1$.\\
(2) $\,a\in {\rm ureg}\,(R)$.\\
(3) $\,a\,$ {\it is a product of elements in} $\,{\rm ureg}\,(R)$.\\
(4) $\;\overline{a}\,$ {\it is unit-regular in
$\,\overline{R}=R/{\rm rad}\,(R)$.}


\bs\nt
{\bf Proof.} We will first prove the equivalence of (1), (2)
and (3). To begin with, $(2)\Rightarrow (3)\,$ is trivial, while
$(3)\Rightarrow (1)$ follows from the Product Theorem in view of
what is said in Example 2.5C. Next, assume (1), and write $\,a=axa\,$
for some $\,x\in R$. Then Theorem 2.2(3) gives an element $\,b\in R\,$
such that $\,a+(1-ax)\,b\in {\rm U}(R)$. Left-multiplying this by
the idempotent $\,e=ax$, we get $\,a\in e\cdot {\rm U}(R)\subseteq
{\rm ureg}\,(R)$, which proves (2).  Finally, $(2)\Rightarrow (4)\,$
is trivial, and $(4)\Rightarrow (1)\,$ follows from Theorem 2.4(D)
together with the already known fact that (4) implies that
$\,{\rm sr}\,(\overline{a})=1\,$ in $\,\overline{R}$. \qed




\bs\nt
{\bf Remark.} In view of Theorem 2.4(A), $(2)\Leftrightarrow (1)$ in
the theorem above shows easily that, for any ring element $\,a\in R\,$
(without any regularity assumption), $\,{\rm sr}\,(a)=1\,$ iff, for any
$\,t\in R$, $\,aR+tR=R\,\Rightarrow\,a+tb\in {\rm ureg}\,(R)\,$ for some
$\,b\in R$.  First stated by Altun and \"{O}zcan in [AO$_1$:~(2.6)],
this fact is an improvement of a result of Zabavsky in [Za:~(7.5)].

\bs\nt
{\bf Example 4.2A.} Over any ring $\,S$, a row vector $\,{\bf v}
=(a_1,\dots,a_n)\in S^n\,$ is called a {\it completable row\/} if
there exists a matrix $\,V\in {\rm GL}_n(S)\,$ with $\,{\bf v}\,$
as its first row. For any such $\,{\bf v}$, the matrix $\,A\in R
={\mathbb M}_n(S)\,$ with first row $\,{\bf v}\,$ and all other rows
zero has stable range one, since $\,A=E_{11}V\in {\rm ureg}\,(R)$.

\mk\nt
{\bf Example 4.2B.} In any matrix ring $\,R={\mathbb M}_2(S)$, let
$\,U=\footnotesize{\begin{pmatrix}1&0\\s&1\end{pmatrix}}\in {\rm U}(R)\,$
and $\,E=\footnotesize{\begin{pmatrix}1&t\\0&0\end{pmatrix}}\in
{\rm idem}\,(R)$. Then $\,B:=UE=\footnotesize{\begin{pmatrix}
1&t\\s&s\,t\end{pmatrix}}$ and $\,A:=EU=\footnotesize{\begin{pmatrix}
1+ts&t\\0&0\end{pmatrix}}$ are both unit-regular in $\,R$, so
$\,{\rm sr}_R(B)={\rm sr}_R(A)=1$.  The last equality here also
follows from Example 4.2A since the first row $\,(1+ts,\,t)\,$
of $\,A\,$ is completable to the invertible matrix
$\,V=\footnotesize{\begin{pmatrix}1+ts&t\\s&1\end{pmatrix}}$. In the
special case where $\,S\,$ is a commutative ring, the fact that
$\,{\rm sr}_R(B)=1\,$ was first noted by C\v{a}lug\v{a}reanu and
Pop in [CP$_1$].

\bs
Recall that an element $\,y\in R\,$ is called a {\it reflexive
inverse\/} of an element $\,a\in R\,$ if $\,a=aya\,$ and $\,y=yay$.
Making use of this standard terminology, we can produce a number of
additional ``inter-connected'' characterizations for unit-regular
elements in any ring $\,R$, as in the following result Theorem 4.3.
Here, the two characterizations (8) and (9) are taken from Wei's
paper [We:~Lemma 2.1], but our proofs below are considerably simpler.

\bs\nt
{\bf Theorem 4.3.} {\it For any ring element $\,a\in R$,
the following statements are equivalent\/}:

\mk\nt
(1) $\;a\in {\rm ureg}\,(R)$.\\
(2) {\it $\;a\in {\rm reg}\,(R)$, and every reflexive inverse $\,x\,$
of $\,a\,$ is unit-regular.}\\
(3) {\it $\;a\in {\rm reg}\,(R)$, and every reflexive inverse $\,x\,$
of $\,a\,$ has $\,{\rm sr}\,(a)=1$.}\\
(4) {\it $\;a\,$ has a reflexive inverse $\,x\in R\,$ such that
$\,{\rm sr}\,(x)=1$.}\\
(5) {\it $\;a\,$ has an inner inverse $\,x\in R\,$ such that
$\,{\rm sr}\,(x)=1$.}\\
(6) {\it $\;a\,$ has a reflexive inverse that is unit-regular.}\\
(7) {\it $\;a\,$ has an inner inverse that is unit-regular.}\\
(8) {\it \;There exist $\,y\in R\,$ and $\,u\in {\rm U}(R)\,$ such
that $\,a=aya=ayu$.}\\
(9) {\it \;There exist $\,y\in R\,$ and $\,w\in {\rm ureg}\,(R)\,$
such that $\,a=aya=ayw$.}\\
(10) {\it There exist $\,y,\,z\in R\,$ with $\,{\rm sr}\,(z)=1\,$
such that $\,a=aya=ayz$.}


\bs\nt
{\bf Proof.} We will first prove the equivalence of (1) through (5).
To begin with, $(2)\Rightarrow (3)\,$ follows from Theorem 4.1,
and $(3)\Rightarrow (4)\Rightarrow (5)\,$ are trivial.

\mk\nt
$(5)\Rightarrow (1)$.  Assume that $\,a=axa\,$ for some $\,x\in R\,$
such that $\,{\rm sr}\,(x)=1$.  According to Theorem 2.2, there
exists $\,b\in R\,$ such that $\,u\!:=x+(1-xa)\,b\in {\rm U}(R)$.
Left multiplying this by $\,a\,$ gives $\,au=ax+(a-axa)\,b=ax$, so
$\,a=(ax)\,a=aua\in {\rm ureg}\,(R)$.

\mk\nt
$(1)\Rightarrow (2)$.  Let $\,x\,$ be any reflexive inverse
of $\,a\in {\rm ureg}\,(R)$. We have $\,{\rm sr}\,(a)=1\,$
by Theorem 4.1, and $\,x\,$ has inner inverse $\,a$. Thus,
applying the preceding implication $(5)\Rightarrow (1)\,$ to
$\,x\,$ shows that $\,x\in {\rm ureg}\,(R)$.

\mk
Having shown the equivalence of (1)--(5), we'll next prove that
$(1)\Rightarrow (6)\Rightarrow (7)\Rightarrow (5)$. Here,
$(6)\Rightarrow (7)$ is a tautology, and $(7)\Rightarrow (5)$
follows from Theorem 4.1.  To prove $(1)\Rightarrow (6)$, fix
a unit inner inverse $\,u\,$ for $\,a\in {\rm ureg}\,(R)$. Then
$\,a\,$ has reflexive inverse $\,uau$, which is easily seen to
be unit-regular.

\mk
To complete the proof of the theorem, we note that $(1)\Rightarrow (8)$
follows by writing $\,a=aya\,$ for some $\,y\in {\rm U}(R)\,$ and
choosing $\,u\,$ to be $\,y^{-1}$. Next, $(8)\Rightarrow (9)\,$ is
trivial, while $(9)\Rightarrow (10)\,$ follows from Theorem 4.1.
Finally, if $\,y\,$ and $\,z\,$ exist as in (10), then $\,a=aya \in
{\rm reg}\,(R)$, and $\,a=(ay)\,z\,$ (with $\,ay\in {\rm idem}\,(R)\,$)
implies that $\,{\rm sr}\,(a)=1\,$ again by Theorem 2.8.  Thus,
$\,a\in {\rm ureg}\,(R)\,$ by Theorem 4.1.\qed




\bs\nt
{\bf Remark 4.4.} With respect to the equivalence of (1) and (6) in the
above theorem, one may wonder what would happen if we try to strengthen
the condition (6) to $\,a\,$ having a reflexive inverse that is
{\it strongly regular.}  This question has been recently answered
by X.~Mary, who showed in [My:~Theorem 4.1] that this strengthened
version of (6) is equivalent to $\,a\,$ having a ``special clean
decomposition'' (see [KNS]) in the sense that $\,a=e+u\,$ where
$\,e=e^2\in R\,,$ and $\,u\in {\rm U}(R)\,$ is such that $\,a=au^{-1}a$.

\bs
Returning to the theme of Theorem 4.1 and using once more the
Product Theorem 2.8, we will next prove the interesting result
below about products of unit-regular elements in any ring $\,R$.

\bs\nt
{\bf Theorem 4.5.} (1) {\it Any finite product of elements
in $\,{\rm ureg}\,(R)\,$ has stable range one.}\\
(2) {\it If $\,a_1,\dots,a_n\in {\rm ureg}\,(R)\,$ are
such that $\,b\!:=a_1\cdots a_n\in {\rm reg}\,(R)$, then
$\,b\in {\rm ureg}\,(R)$.}  ({\it In particular, if
$\,a\in {\rm ureg}\,(R)\,$, then for any positive integer $\,n$,
$\,a^n\in {\rm reg}\,(R)\,$ iff $\,a^n\in {\rm ureg}\,(R)$.})\\
(3) {\it If the set $\,{\rm reg}\,(R)\,$ is closed under
multiplication, then so is the set $\,{\rm ureg}\,(R)$.}

\bs\nt
{\bf Proof.} (1) follows from Theorem 4.1 and the Product
Theorem 2.8. Since obviously (2) implies (3), we may finish by
proving (2). By (1), we have $\,{\rm sr}\,(b)=1$. If we assume
that $\,b\in {\rm reg}\,(R)$, another application of Theorem 4.1
shows that $\,b\in {\rm ureg}\,(R)$. (Obviously, part (3)
here implies the main result of Hannah and O'Meara [HO:~(3.2)],
which states that, in any {\it regular ring\/} $\,R$, the set
$\,{\rm ureg}\,(R)\,$ is closed under multiplication.  For an
independent proof of (2) in the special case where $\,n=2$,
see [La$_{3}$: (6.33)].)  \qed


\bs\nt
{\bf Remark 4.6.} Concerning Theorem 4.5(2) again, we
should point out that, if $\,a_1,\dots,a_n\in {\rm sreg}\,(R)\,$
are such that $\,b\!:=a_1 \cdots a_n\in {\rm reg}\,(R)$, then
$\,b\,$ may not be in $\,{\rm sreg}\,(R)\,$ (although we do have
$\,b\in {\rm ureg}\,(R)\,$ by Theorem 4.5(1)). For instance, in
the matrix ring $\,R={\mathbb M}_2\bigl({\mathbb Z}\bigr)$, both
$\,a_1=\footnotesize{\begin{pmatrix}1&2\\0&0\end{pmatrix}}$
and $\,a_2=\footnotesize{\begin{pmatrix}0&-1\\1&1\end{pmatrix}}$
are in $\,{\rm sreg}\,(R)$, but
$\,b\!:=a_1a_2=\footnotesize{\begin{pmatrix}2&1\\0&0\end{pmatrix}}
\in {\rm ureg}\,(R)\,$ is not in $\,{\rm sreg}\,(R)\,$ since
$\,b^2=2\,b\,$ clearly implies that $\,b\notin b^2 R$.

\bs
To state a consequence of Theorem 4.5, we recall from [KL$_2$] that
a ring $\,R\,$ is said to be IC (``internally cancellable'') if,
whenever $\,R_R=A\oplus B=A'\oplus B'$, $\,A\cong A'\,\Rightarrow\,
B\cong B'\,$ in the category of right $\,R$-modules.  From
[KL$_2$:~(1.5)], it is well known that this property is left-right
symmetric, and is equivalent to the condition that $\,{\rm reg}\,(R)
={\rm ureg}\,(R)$. In view of this ``element-wise'' characterization
of IC rings, the following result is an immediate consequence of
Theorem 4.1 and Theorem 4.5(1).

\bs\nt
{\bf Corollary 4.7.} {\it A ring $\,R\,$ is {\rm IC} iff every
element $\,a\in {\rm reg}\,(R)\,$ has stable range one,  iff
every element in $\,{\rm reg}\,(R)\,$ is a product of elements
in $\,{\rm ureg}\,(R)$.}


\bs\nt
{\bf Remark 4.8.} If $\,a,b\in {\rm ureg}\,(R)$, there is also
a known ideal-theoretic criterion for $\,ab\,$ to be regular
(or unit-regular) in $\,R\/$; namely, that the right ideal
$\,{\rm ann}_r(a)+bR\,$ be a direct summand of $\,R_R$. For
a proof of this interesting criterion, see [La$_4$:~Theorem 6.33].
Having proved Corollary 4.7, we should mention that the rings $\,R\,$
with the property that $\,{\rm reg}\,(R)\,$ is closed under
multiplication are called {\bf reg-closed rings\/} in [La$_4$:~\S6],
where it is shown that these are precisely the {\bf SSP rings\/}
introduced earlier by Garcia in [Ga]. (The term ``SSP'' was in
reference to a certain ``summand sum property'' characterizing
such rings.) For a fuller discussion on these rings and a more
complete listing of the relevant literature on them, see again
[La$_4$:~\S6].


\bs
For any commutative ring $\,S\,$ with only trivial idempotents,
it is known from [KL$_2$:~Example 5.9(1)] that the matrix ring
$\,{\mathbb M}_2(S)\,$ is IC. In fact, this conclusion holds even
without the assumption that $\,{\rm idem}\,(S)=\{0,1\}$, since
a recent argument of Swan reported in [La$_4$:~Theorem 23.16]
showed that $\,{\mathbb M}_2(S)\,$ is IC as long as $\,S\,$ is
a commutative ring.  Assuming this nice result, we can deduce
the following.

\bs\nt
{\bf Corollary 4.9.} {\it For any commutative ring $\,S$,
any (von Neumann) regular matrix in $\,{\mathbb M}_2(S)\,$
is unit-regular, and hence has stable range one.}

\bs\nt
{\bf Example 4.10.} Even for a commutative {\it domain\/} $S$,
${\mathbb M}_3(S)$ may fail to be IC.  For instance,
it is pointed out after the proof of [La$_4$:~(23.16)] that, for
$\,S={\mathbb R}\,[x,y,z]/(x^2+y^2+z^2-1)\,$ (the coordinate ring
of the real $2$-sphere), the matrix ring $\,R={\mathbb M}_3(S)\,$
is {\it not\/} IC, with the matrix $\,A=\footnotesize{\begin{pmatrix}
x&y&z\\0&0&0\\0&0&0\end{pmatrix}}$ being regular (with reflexive
inverse $\,A^T\,$) but not unit-regular, so that $\,{\rm sr}_R(A)
\neq 1$. Finally, noting that $\,A\,$ has a factorization $\,BC\,$
where $\,B=\footnotesize{\begin{pmatrix}1&y&z\\0&0&0\\0&0&0
\end{pmatrix}}$ is an idempotent and $\,C={\rm diag}\,(x,1,1)$,
we see that $\,{\rm sr}_R(C)\neq 1\,$ by Theorem 2.8.  This is
in sharp contrast to our forthcoming result Theorem 7.1, which
implies that any diagonal matrix $\,D={\rm diag}\,(a,b,0)\,$ (over
any base ring $\,S\,$) has stable range one in $\,{\mathbb M}_3(S)$.


\bs
To close this section, we shall prove one more result relating the
notion of IC (``internal cancellation'') to the ring-theoretic notion
of stable range one. {\it In the case of an exchange ring $\,R$,} it has
been well known since 1995 that $\,{\rm sr}\,(R)=1\,$ iff $\,R\,$ is IC,
iff $\,{\rm reg}\,(R)={\rm ureg}\,(R)$.  This nice result was first
proved by Yu in [Yu$_1$:~Theorem 9] and Camillo-Yu in [CY$_1$:~Theorem 3],
and later recapitulated by Chen in [Ch$_2$, Ch$_3$], and Khurana-Lam
in [KL$_2$].  While this classical result was purely expressed in terms
of the ``global'' ring-theoretic properties of $\,R$, the element-wise
approach to stable range one taken in our present paper would seem to
suggest that there should exist a corresponding criterion characterizing
elements $\,a\in R\,$ with $\,{\rm sr}\,(a)=1\,$ in the case where,
say, all elements in $\,aR\,$ are {\it suitable\/} (in the sense of
Nicholson [Ni$_1$]).  Optimally, such an element-wise criterion should
be expressible in terms of ``certain'' von Neumann regular elements
{\it associated with $\,a\in R\,$} being automatically unit-regular.
Moreover, it should naturally imply the classical Camillo-Yu theorem
in the case where $\,R\,$ itself is an exchange ring. Luckily, such an
element-wise criterion does exist, as we will show in the result below.

\bs\nt
{\bf Theorem 4.11.} {\it Let $\,a\in R\,$ be such that all elements
in $\,aR\,$ are suitable. Then $\,{\rm sr}\,(a)=1\,$ iff, for any
$\,f\in {\rm idem}\,(R)$, $\,fa\in {\rm reg}\,(R)\Rightarrow
fa\in {\rm ureg}\,(R)$.}

\bs\nt
{\bf Proof.} First assume $\,{\rm sr}\,(a)=1$. If $\,f\in
{\rm idem}\,(R)\,$ is such that $\,fa\in {\rm reg}\,(R)$, then
$\,{\rm sr}\,(fa)=1\,$ by Theorem 2.8, so Theorem 4.1 implies
that $\,fa\in {\rm ureg}\,(R)$.  Conversely, assume that for any
$\,f\in {\rm idem}\,(R)$, $\,fa\in {\rm reg}\,(R)\Rightarrow fa\in
{\rm ureg}\,(R)$. To show that $\,{\rm sr}\,(a)=1$, we'll verify
the criterion (7) in Corollary 2.3; that is, given any equation
$\,ax+e=1\,$ with $\,e=e^2$, we want to find some $\,b\in R\,$ such
that $\,a+e\,b\in {\rm U}(R)$.  Letting $\,f\!:=1-e=ax$, we have
$\,(fa)\,x\,(fa)=f^3 a=fa$, so $\,fa\in {\rm reg}\,(R)$. By assumption,
this implies that $\,fa\in {\rm ureg}\,(R)$, so $\,{\rm sr}\,(fa)=1$.
From $\,(fa)\,x+e=f^2+e=1$, it then follows that $\,fa+er\in
{\rm U}(R)\,$ for some $\,r\in R$. This amounts to $\,a+e\,(r-a)
\in {\rm U}(R)$, as desired.\qed

\mk

\bs\nt
{\bf \S5. \ Stable Range One Elements Via Nilpotency}

\bs
Since units and idempotents in any ring $\,R\,$ have stable range
one, it is natural to wonder which nilpotent elements $\,a\in
{\rm nil}\,(R)\,$ might have the same property. To give some
examples of such elements, let us first recall the following two
classical element-wise notions that were closely related to the
notion of nilpotent elements in rings. First, an element $\,a\in R\,$
is called {\bf strongly nilpotent\/} (after J.~Levitzki [Lv]) if
every sequence $\,a_1,a_2,a_3,\dots\,$ such that $\,a_1=a\,$
and $\,a_{n+1}\in a_nRa_n\,$ for all $\,n\,$ is eventually zero.
Secondly, an element $\,a\in R\,$ is called {\bf quasi-nilpotent\/}
(after R.~Harte [Ha]) if $\,1-as\in {\rm U}(R)\,$ for every $\,s\in R\,$
commuting with $\,a$.  With these two standard definitions in place,
our first main result in this section can be stated as follows.

\bs\nt
{\bf Theorem 5.1.} {\it Suppose an element $\,a\in R\,$ satisfies
one of the following conditions\/}:

\mk\nt
(A) {\it $\,a\,$ is strongly nilpotent\,};\\
(B) {\it $\,a\,$ is central and quasi-nilpotent\,};\\
(C) {\it $\,a\in {\rm nil}\,(R)$, all elements in $\,aR\,$ are
suitable, and $\,ae=ea\,$ for every $\,e\in {\rm idem}\,(R)$.}

\mk\nt
{\it Then $\,a\in {\rm rad}\,(R)$, and hence $\,{\rm sr}\,(a)=1\,$
by\/} Theorem 2.6(2).

\bs\nt
{\bf Proof.} (A) If $\,a\in R\,$ is strongly nilpotent, Levitzki's result
[Lv:~Theorem 1] guarantees that $\,a\,$ is in the lower nil-radical
of $\,R$.  (For an easily accessible proof for this fact, see [La$_2$:
Exercise 10.17].)  In particular, we have $\,a\in {\rm rad}\,(R)$.

\mk\nt
(B) If $\,a\,$ is central and quasi-nilpotent, then $\,1-as\in
{\rm U}(R)\,$ for every $\,s\in R\,$ (since we do have $\,as=sa\,$
here by assumption).  By [La$_1$:~Lemma 4.1], this shows that
$\,a\in {\rm rad}\,(R)$.

\mk\nt
(C) We first show that $\,aR\cap {\rm idem}\,(R)=\{0\}$.  Indeed,
if $\,e=ax\in {\rm idem}\,(R)\,$ for some $\,x\in R$, then an easy
induction shows that $\,e=a^n x^n\,$ (noting that $\,e=a^n x^n
\Rightarrow e=e(ax)=aex=a(a^n x^n)x=a^{n+1}x^{n+1}$). Thus, taking
$\,n\,$ to be a large integer shows that $\,e=0$. To prove that
$\,a\in {\rm rad}\,(R)$, it suffices to show that $\,(1-ay)\,R=R\,$
for every $\,y\in R\,$ (again by [La$_1$: Lemma 1.4]). Since $\,ay\,$
is a suitable element (by assumption), there exists an idempotent
$\,f\in ayR\,$ such that $\,1-f\in (1-ay)R$. By the first statement
in this paragraph, we must have $\,f=0$, and so $\,1\in (1-ay)R$,
as desired. \qed



\bs
While the elements $\,a\in R\,$ studied in Theorem 5.1 all turned
out to be in $\,{\rm rad}\,(R)$, we will next study the case
of {\it regular nilpotent elements\,}; that is, elements $\,a\in
{\rm reg}\,(R)\cap {\rm nil}\,(R)$.  This time, $\,a\,$ is
no longer in $\,{\rm rad}\,(R)\,$ unless $\,a=0\,$ (according
to [La$_4$:~Theorem 2.20]).  For such elements, we have the
following two known results (1) and (2) from the literature.

\bs\nt
{\bf Theorem 5.2.} {\it Let $\,a\in {\rm reg}\,(R)$, with
$\,a^n=0\,$ for a given integer $\,n\geq 1$.}

\mk\nt
(1) {\it If $\,a^i\in {\rm reg}\,(R)\,$ for all $\,i\geq 1$, then
$\,a\,$ is a product of $\,n\,$ idempotents in $\,R\,$}. {\it In
this case, $\,a^i\in {\rm ureg}\,(R)\,$ and $\,{\rm sr}\,(a^i)=1\,$
for all $\,i\geq 1$.}\\
(2) {\it If $\,R\,$ is an exchange ring, then $\,a\in {\rm ureg}\,(R)$,
and $\,{\rm sr}\,(a^i)=1\,$ for all $\,i\geq 1$.}


\bs\nt
{\bf Proof.} (1) is essentially a result of Beidar, Hannah, O'Meara
and Raphael; see [HO:~(1.2)] and [BOR:~(3.7)].  For an exposition
on this result, see [La$_4$:~(2.46), (6.40)].  Here, we will offer
a quick and completely self-contained proof for (1) in the special
case where $\,n=2$.  Writing $\,a=axa\,$ for some $\,x\in R$, let
$\,e=ax\in {\rm idem}\,(R)$, and $\,f=1-e$. Then $\,ae=a\,(ax)=0$, so
$\,(e+a)\,f=a\,(1-e)=a$. Here, $\,(e+a)^2=e+ea=e+a$, so $\,a\,$ is a
product of two idempotents. Also, $\,a\,(1+fxe)\,a=a\,(1-e)\,xa=a\,$
shows that $\,a\in {\rm ureg}\,(R)\,$ since $\,1+fxe\,$ is a
unipotent unit.  Thus, $\,{\rm sr}\,(a)=1\,$ by Theorem 4.1.
Finally, in the case (2) where $\,R\,$ is assumed to be an exchange
ring, the fact that $\,a\in {\rm reg}\,(R)\cap {\rm nil}\,(R)
\Rightarrow a\in {\rm ureg}\,(R)\,$ is a significant result of Ara
from [Ar:~Theorem 2]. In this case, it follows from Theorem 4.1 and
Theorem 2.8 that $\,{\rm sr}\,(a^i)=1\,$ for all $\,i\geq 1$.\qed



\bs
In practice, there are many concrete examples of non-central nilpotent
elements $\,a\in R\,$ with the property that $\,{\rm sr}_R(a)=1$.  In
Example 5.3 below, we will work in a $\,2\times 2\,$ matrix ring to
produce some such elements that are even unit-regular.

\bs\nt
{\bf Example 5.3.} For any square-zero element $\,s\in S$, consider
the square-zero matrix $\,A\!:=\footnotesize{\begin{pmatrix}s&1\\0&-s
\end{pmatrix}} \in R={\mathbb M}_2(S)$. This matrix has a unit inner
inverse $\footnotesize{\begin{pmatrix}1&0\\1&1 \end{pmatrix}}$.
Thus, $\,A\in {\rm ureg}\,(R)$, and so $\,{\rm sr}_R(A)=1\,$
by Theorem 4.1.  More generally, for any $\,t\in S$, the
``equivalent'' matrix $\,B\!:=A\,\footnotesize{\begin{pmatrix}
1&t\\0&1\end{pmatrix}}=\footnotesize{\begin{pmatrix}s&1+st\\0&-s
\end{pmatrix}}$ continues to be unit-regular.  So again
$\,{\rm sr}_R(B)=1$, and hence $\,{\rm sr}_R(B^2)=1\,$ too. Here,
$\,B^2=-sts\,E_{12}\,$ may not be zero, although we do have $\,B^3=0$.



\bs
In general, if $\,a\,$ is just a nilpotent element in a ring $\,R$,
we may not have $\,{\rm sr}_R(a)=1\,$ even in the case where $\,a\in
{\rm reg}\,(R)$. Indeed, in [NS:~Theorem 3.19] and [ArO:~Theorem 2.1],
an example of a ring element $\,a\in {\rm reg}\,(R)\setminus
{\rm ureg}\,(R)\,$ (in some ring $\,R\,$) was given such that
$\,a^3=0$, so in particular $\,{\rm sr}_R(a)\neq 1\,$ by Theorem 4.1.
Without assuming regularity, on the other hand, we will offer below
an example of a square-zero element $\,a\,$ in an abelian ring
$\,R\,$ such that $\,{\rm sr}_R(a)\neq 1$.

\bs\nt
{\bf Example 5.4.} Let $\,R\,$ be the ring $\,k\,\langle a,x\rangle\,$
generated over the $2$-element field $\,k={\mathbb F}_2\,$ with a
single relation $\,a^2=0$. A free $\,k$-basis of $\,R\,$ is given by
the set of (noncommutative) monomials in the ``variables'' $\,a\,$ and
$\,x\,$ which do not involve a factor $\,a^2$. Thus, $R\,$ has an obvious
$\,{\mathbb Z}^{+}$-graded ring structure. From this, we see easily
that $\,{\rm idem}\,(R)=\{0,1\}$; in particular, $R\,$ is an abelian
ring. Assuming a result of Bergman reported in [DL:~(2.5)], the unit
group $\,{\rm U}(R)\,$ is given by $\,1+k\,a+aR\,a\,$. (In particular,
all units in $\,R\,$ are unipotents, so $\,R\,$ is a ``UU ring''
in the sense of C\v{a}lug\v{a}reanu [Ca].)  {\it We claim that
$\,{\rm sr}\,(a)\neq 1$.}  Indeed, if otherwise, Theorem 2.2 would
imply the existence of an element $\,b\in R\,$ such that $\,a+b+axb
\in {\rm U}(R)$.  Writing out $\,b\,$ as a ``noncommutative polynomial''
in $\,a\,$ and $\,x\,$ in which $\,a^2\,$ does not appear, we see
that $\,b\,$ must have constant term $\,1$.  If $\,n\geq 0\,$ is
chosen largest such that $\,b\,$ has a monomial term $\,(ax)^n$,
then $\,a+b+axb\,$ would have a monomial term $\,(ax)^{n+1}$. This
would contradict the fact that $\,{\rm U}(R)=1+k\,a+aR\,a$, so we
have proved that $\,{\rm sr}\,(a)\neq 1$.

\bs
Going back to our earlier work (in the proof of Theorem 5.2(1))
on a square-zero regular element $\,a=axa\,$ in any ring $\,R$,
recall that, for the idempotent $\,e=ax\,$ and its complementary
idempotent $\,f=1-e\,$ there, we have $\,ae=0\,$ and so
$\,a=af=eaf\in eRf$. The conclusion that $\,{\rm sr}\,(a)=1\,$
in that situation prompted us to the discovery of the following
more general result.

\bs\nt
{\bf Theorem 5.5.} {\it Let $\,e,f\in {\rm idem}\,(R)\,$ be such
that $\,ef=fe=0$. Then $\,{\rm sr}\,(erf)=1\,$ for all $\,r\in R$,
and $\,{\rm sr}\,(erfse)=1\,$ for all $\,r,s\in R$.}

\bs\nt
{\bf Proof.} It is easy to check that $\,f+erf\in {\rm idem}\,(R)$.
This shows that $\,erf=e\,(f+erf)\,$ is a product of two idempotents
in $\,R$. Since idempotents have stable range one, the Product
Theorem 2.8 implies that $\,{\rm sr}\,(erf)=1$. We have similarly
$\,{\rm sr}\,(fse)=1$, so another application of the Product
Theorem gives $\,{\rm sr}\,(erfse)=1$.\qed


\bs\nt
{\bf Corollary 5.6.} {\it Let $\,p,q\in R\,$ be such
that $\,pq=qp=0\,$ and $\,u\!:=p+q\in {\rm U}(R)$. Then
for $\,a\!:=prq\,$ where $\,r\,$ is any element of $\,R$,
we have $\,{\rm sr}\,(a)=1$.}

\bs\nt
{\bf Proof.} Note that the elements $\,p,q,u\,$ pairwise commute.
If we define $\,e=u^{-1}p\,$ and $\,f=u^{-1}q$, then $\,ef=fe=0\,$
and $\,e+f=1$, so $\,e,f\,$ are complementary idempotents in $\,R$.
Since $\,p=ue\,$ and $\,q=uf$, we have $\,u^{-1}a=u^{-1}prq=eruf$.
Thus, Theorem 5.5 gives $\,{\rm sr}\,(u^{-1}a)=1$, and hence
$\,{\rm sr}\,(a)=1\,$ by Theorem 2.4(B).\qed

\bs
In some other special situations, we may also exploit the conclusion
in Theorem 5.5 to get useful information on $\,{\rm sr}_R(r)\,$
for an element $\,r\,$ in a corner ring $\,eR\,e$.  To do this, let
us first recall the following basic notions in general ring theory.
First, two elements $\,x,y\in R\,$ are said to be {\it equivalent\/}
if $\,y=uxv\,$ for some $\,u,v\in {\rm U}(R)$. Second, two idempotents
$\,e,g\in R\,$ are said to be {\it isomorphic\/} if $\,R\,e\cong Rg\,$
as left $\,R$-modules; or equivalently, if $\,R\,e\cong Rg\,$ as right
$\,R$-modules. By standard ring theory (see, e.g.~[La$_1$:~Prop.\,21.20]),
this amounts to the existence of two elements $\,a\in eRg\,$ and
$\,b\in gRe\,$ such that $\,e=ab\,$ and $\,g=b\,a$. In terms of the
two basic notions above, our next result gives a {\it sufficient\/}
condition for every element in a corner ring $\,eR\,e\,$ to have
stable range one in a given ring $\,R$.

\bs\nt
{\bf Theorem 5.7.} {\it Let $\,e,\,f\,$ be a pair of complementary
idempotents in $\,R\,$ such that $\,e\,$ is isomorphic to some
idempotent $\,g\in fRf$.  Then any element $\,r\in eR\,e\,$ is
equivalent to some element in $\,eRf$.  In particular,
$\,{\rm sr}_R(r)=1$.}

\bs\nt
{\bf Proof.} Given that $\,e\,$ and $\,g\,$ are isomorphic
idempotents, we may choose the two elements $\,a,b\in R\,$
as in the paragraph preceding the statement of this theorem.
Let $\,h=f-g\in {\rm idem}\,(fRf)$. Upon expanding $\,(a+b+h)^2\,$
into nine terms, six of these terms are zero, leaving us with
$$
(a+b+h)^2=ab+ba+h^2=e+g+h=1.
$$
In particular, $\,a+b+h\in {\rm U}(R)$. Thus, $\,r\in eR\,e\,$
is equivalent to
$$
r\,(a+b+h)=ra\in (eR\,e)\,(eRg)\subseteq eRf,
$$
as claimed.  Given this information, the last conclusion that
$\,{\rm sr}_R(r)=1\,$ follows from Theorem 5.5 and Theorem 2.4(B).\qed


\bs
Note that the theorem above may no longer hold in general if we
try to remove the hypothesis on $\,e\,$ being isomorphic to some
idempotent $\,g\in fRf$. For instance, let $\,R\,$ be the commutative
ring $\,{\mathbb Z}\times {\mathbb Z}$, and let $\,e=(1,0)$,
$\,f=(0,1)=1-e$. Here, $\,e\,$ is certainly not isomorphic to any
idempotent in $\,fRf$. For any integer $\,n\geq 2$, $\,(n,0)\in eR\,e
\cong {\mathbb Z}\,$ does not have stable range one in $\,R$, for
otherwise the projection map from $\,R\,$ to $\,eR\,e\,$ would have
given $\,{\rm sr}_{\mathbb Z}(n)=1$, which is not the case according
to Example 2.5B.

\bs
Applying Theorem 5.5 to the case of a full matrix ring
$\,R={\mathbb M}_n(S)$, we will next prove the following very
useful result which produces many matrices of stable range one
in such a matrix ring. Here, as usual, the $\,E_{ij}$'s denote
the matrix units in $\,R$.

\bs\nt
{\bf Theorem 5.8.} {\it Let $\,A=(a_{ij})\in R={\mathbb M}_n(S)\,$
where $\,n\geq 2$, and assume that, for some integer $\,k\in [1,n]$,
$\,a_{ij}=0\,$ for all $\,i\neq k\,$ and all $\,j$, and that $\,a_{kj}
=0\,$ for some $\,j$. Then $\,{\rm sr}\,(A)=1$.  In particular,
$\,{\rm sr}\,\bigl(aE_{kj}\bigr)=1\,$ for all $\,a\in S\,$ and all
$\,k,j$.}

\bs\nt
{\bf Proof.} After a row permutation, we may assume in view of
Theorem 2.4(B) that $\,k=1$. Following this by a column permutation,
we may likewise assume that $\,a_{11}=0$. For $\,e\!:=E_{11}$ and
$\,f=I_n-e$, we have then $\,A\in eRf$. Therefore, Theorem 5.5
implies that $\,{\rm sr}\,(A)=1$.\footnote{In 2020, after an
initial announcement of our results, a somewhat different
proof for $\,{\rm sr}\,\bigl(aE_{kj}\bigr)=1\,$ was given
by C\v{a}lug\v{a}reanu and Pop in [CP$_1$:~Proposition 7].} \qed


\bs\nt
{\bf Remark 5.9.} Of course, the last statement in Theorem 5.8 implies
that, as long as $\,n\geq 2$, every matrix in $\,{\mathbb M}_n(S)\,$
is a sum of $\,n^2\,$ matrices of stable range one.  We did not state
this as a part of Theorem 5.8, since there is already a much better
result in the literature.  According to Henriksen [He$_2$], every
matrix in $\,{\mathbb M}_n(S)\,$ is in fact a sum of {\it three\/}
invertible matrices, each of which has stable range one by Example 2.1.
The foregoing discussion brings into focus the following question:
which rings $\,R\,$ could have the property that $\,R\,$ is additively
spanned by the set of stable range one elements in $\,R\,$?  Examples
of such rings include: (1) all rings of stable range one (for instance
those mentioned earlier in Example 2.5A); (2) the ring $\,{\mathbb Z}\,$;
and (3) the matrix rings $\,{\mathbb M}_n(S)\,$ for $\,n\geq 2$). Other
examples include (4) all clean rings in the sense of Nicholson [Ni$_1$],
and more generally, (5) any ring in which every element is a sum of
units and idempotents (e.g.~${\mathbb Z}$, $\,{\mathbb Z}\,[\sqrt{2}\,]$,
$\,{\mathbb Z}\,[\sqrt{3}\,]$, or any group ring $\,k\,[G]\,$
where $\,k\,$ is a ring in which every element is a sum of units).
On the other hand, any polynomial ring $\,R={\mathbb Z}\,[x_1,\dots,x_n]\;
(n\geq 1)\,$ is easily seen to be a ``non-example'', since a quick
application of Theorem 2.4(D) shows that the set of stable range one
elements in $\,R\,$ is $\,\{0,\,{\pm 1}\}$, whose additive span is
$\,{\mathbb Z}\subsetneq R$.


\bs\nt
{\bf Example 5.10.} Using Theorem 5.8, it is easy to come up
with examples of ring elements of stable range one that are
not unit-regular (or even regular).  For instance, in the ring
$\,R={\mathbb M}_2\bigl({\mathbb Z}\bigr)$, $\,A=2E_{11}\,$ has
stable range one by Theorem 5.8, but $\,A\notin {\rm reg}\,(R)\,$
since any equation $\,A=ABA \in R\,$ for $\,B=(b_{ij})\in R\,$
would lead to a ``bad'' equation $\,2=4\,b_{11}\in {\mathbb Z}$.


\bs\nt
{\bf Example 5.11.} This example is a generalization of a result
of C\v{a}lug\v{a}reanu and Pop [CP$_1$:~Theorem 9] on $\,2\times 2\,$
matrices over B\'{e}zout rings. Let $\,R={\mathbb M}_2(S)\,$ where
$\,S\,$ is a commutative ring, and let $\,p,q\in S\,$ be such that
$\,pS+qS=aS\,$ for some $\,a\in S\,$ that is not a $\,0$-divisor
in $\,S$. We claim that the matrix $\,C=\footnotesize{\begin{pmatrix}
p&q\\0&0\end{pmatrix}}\in R\,$ has stable range one.  To see this,
we write $\,p=as,\;q=at\,$ for some $\,s,t\in S$, and $\,a=px-qy\,$
for some $\,x,y\in S$.  Since $\,a\,$ is not a $\,0$-divisor,
we have $\,sx-ty=1$, so $\,U\!:=\footnotesize{\begin{pmatrix}
s&t\\y&x\end{pmatrix}}\in {\rm U}(R)$. From the factorization
$\,C=(aE_{11})\,U$, we see that $\,{\rm sr}_R(C)=1\,$ by using
Theorem 2.4(B) and Theorem 5.8.

\bs\nt
{\bf Example 5.12.} The conclusion that $\,{\rm sr}\,\bigl(aE_{ij}\bigr)
=1\,$ in Theorem 5.8 can be more useful than we might have initially
thought. For instance, applying it to a square-zero matrix
$$
\mbox{$M=\footnotesize{\begin{pmatrix}0&0&a&b\\0&0&c&d\\0&0&0&0\\0&0&0&0
  \end{pmatrix}}$, viewed as a block matrix of the form
  $\footnotesize{\begin{pmatrix}{\bf 0}_2&N\\{\bf 0}_2 & {\bf 0}_2
  \end{pmatrix}}$}
$$
where $\,{\bf 0}_2\,$ denotes the $\,2\times 2\,$ zero matrix
and $\,N\in {\mathbb M}_2(S)$, we may conclude from
Theorem 5.8 that $\,{\rm sr}\,(M)=1\,$ in any matrix ring
$\,{\mathbb M}_4(S)={\mathbb M}_2\bigl({\mathbb M}_2(S)\bigr)$.
Similarly, we see that any block diagonal matrix
$\,{\rm diag}\,(N,\,{\bf 0}_2)\,$ has stable range one
in $\,{\mathbb M}_4(S)$. However, Theorem 5.8 does not give any
information about the stable range of, say, $\,{\rm diag}\,(x,z,0)\,$
in $\,{\mathbb M}_3(S)$. To handle this case, we will have to use
another technique; see Example 7.1 in \S7.

\bs
While in this section we have mainly focused our attention on the
stable range one problem for nilpotent elements, we should make the
following supplemental remark about some related results in the
literature. Using an element-wise version of an argument developed
by Goodearl and Menal for proving their ``Theorem 5.8'' in [GM],
one can show that, if $\,a\in R\,$ is {\bf strongly $\,\pi$-regular\/}
(in the sense that $\,a^n\in Ra^{n+1}\cap a^{n+1}R\,$ for some
$\,n\geq 1$), and if additionally $\,a^n\in {\rm reg}\,(R)\,$ for
all $\,n\geq 1$, then $\,a\in {\rm ureg}\,(R)\,$ (and hence
$\,{\rm sr}\,(a)=1\,$); see [CY$_1$:~Theorem 5], [BOR: Corollary 3.7],
[Khu: Theorem 4], or [La$_4$:~Theorem 5.24].  However, the condition
$\,a^2 \in {\rm reg}\,(R)\,$ may already fail (even in the case
where $\,a\in {\rm ureg}\,(R)\,$ with $\,a^3=0\,$), as was shown
by  Camillo-Yu
in [CY$_1$:~Example 9], Yu in [Yu$_2$], and Patr\'icio-Hartwig in [PaH].  
For a self-contained exposition on this,
see [La$_4$:~Example 1.12C].


\mk

\bs\nt
{\bf \S6. Stable Range One Elements Via Peirce Decompositions}

\bs
We offer in this section some more constructive results on stable
range one elements in the context of Peirce decompositions with
respect to corner rings.  To prepare ourselves for the Suspension
Theorem 6.2, we start by first proving the following simple lemma.

\bs\nt
{\bf Lemma 6.1.} {\it Let $\,e,f\in {\rm idem}\,(R)\,$ be such
that $\,e+f=1$, and let $\,x\in {\rm U}(eRe)$, $\,y\in {\rm U}(fRf)$.
Then $\,x+p+y\in {\rm U}(R)\,$ for every $\,p\in fR\,e$.}

\bs\nt
{\bf Proof.} Let $\,x'\,$ be the inverse of $\,x\,$ in $\,eRe$,
and let $\,y'\,$ be the inverse of $\,y\,$ in $\,fRf$. Letting
$\,p'=-y'px'\in fR\,e$, we have $\,yp'=-yy'px'=-fpx'=-px'$, so
$$
(x+p+y)\,(x'+p'+y')=x\,x'+p\,x'+y\,p'+y\,y'=x\,x'+y\,y'=e+f=1,
$$
and similarly, $\,(x'+p'+y')\,(x+p+y)=1$.  This shows that
$\,x'+p'+y'\,$ is an inverse of $\,x+p+y\,$ in $\,R$.  (In the
more suggestive notation of Peirce decompositions, the ``Peirce
matrix'' $\,\footnotesize{\begin{pmatrix}x&0\\p&y\end{pmatrix}}$
is invertible, with inverse given by the Peirce matrix
$\,\footnotesize{\begin{pmatrix}x'&0\\p'&y'\end{pmatrix}}$ since
$\,px'+yp'=0=p'x+y'p$.)  \qed


\bs
The first main result of this section is the ``Suspension Theorem''
below for ring elements of stable range one in a corner ring context.
In the ``iff'' statement of this theorem, the ``if'' part was
essentially first discovered (in the special case $\,p=0\,$) by
Vaserstein in [Va$_2$:~Theorem 2.8].  Although Vaserstein did not
stress (or even talk about) the stable range one notion {\it for
elements,} a careful examination of the proof for his ``Theorem 2.8''
will reveal that he had actually proved the ``if'' part of the
result below {\it in the special case where $\,p=0$.}  Here, we
offer a new proof for the ``if'' part (for all $\,p\in fR\,e$)
using rather different ideas. Also, we will complete this ``if''
part by proving its converse, which then makes for a rather
satisfactory ``if and only if'' statement in the theorem below.
The name of this theorem came about since we may think of the
process of adding the idempotent $\,p+f\,$ to $\,a\in eR\,e\,$
as a kind of ``suspension'' of the element $\,a$.

\bs\nt
{\bf Suspension Theorem 6.2.} {\it Let $\,e,f\in {\rm idem}\,(R)\,$
be such that $\,e+f=1$, and let $\,a\in eR\,e$, $\,p\in fR\,e$.
Then $\,{\rm sr}_{eR\,e}(a)=1\,$ iff $\,{\rm sr}_R(a+p+f)=1$.
In this case, we have necessarily $\,{\rm sr}_R(a)=1$.}

\bs\nt
{\bf Proof.} It suffices to prove the ``iff'' statement since,
with $\,p=0$, $\,{\rm sr}_R(a+f)=1\,$ and $\,{\rm sr}_R(e)=1\,$
 imply that $\,1={\rm sr}_R(e\,(a+f))={\rm sr}_R(a)\,$ by
the Product Theorem 2.8.

\mk
First, assume that $\,{\rm sr}_{eR\,e}(a)=1\,$ and $\,p\in fR\,e$.
To show that $\,{\rm sr}_R(a+p+f)=1$, we take any $\,s\in R$, and
want to produce an element $\,r\in R\,$ satisfying
$$
a+p+f+(1-(a+p+f)\,s)\,r \in {\rm U}(R). \leqno (6.3)
$$
Take $\,r\in eR\,e\,$ such that $\,a+(e-a\,(ese))\,r
\in {\rm U}(eRe)$. We have
\begin{eqnarray*}
  a+p+f+(1-(a+p+f)\,s)\,r &=& a+p+f+r-asr-psr-fsr \\
                    &=& (a+er-aeser)+(p-psr-fsr)+f.
\end{eqnarray*}
Here, $\,a+er-aeser \in {\rm U}(eRe)$, $p-psr-fsr\in fR\,e$,
and $\,f\in {\rm U}(fRf)$. Applying Lemma 6.1, we see that
$$
a+p+f+(1-(a+p+f)\,s)\,r\in {\rm U}(R).
$$
This shows that $\,{\rm sr}_R(a+p+f)=1$.

\mk
Conversely, assume that
$\,{\rm sr}_R(a+p+f)=1$. Since $\,1-p\in {\rm U}(R)$, we have
$\,1={\rm sr}_R((a+p+f)\,(1-p))={\rm sr}_R(a+f)$. For any
$\,s\in eR\,e$, there exists $\,t\in R\,$ such that
$\,\alpha\!:=a+f+(1-(a+f)\,(s+f))\,t \in {\rm U}(R)$.
By direct computation, we have
\begin{eqnarray*}
  \alpha &=& a+f+(e-as)\,t=a+f+(e-as)\,(ete+etf+fte+ftf)  \\
  &=& [\,a+(e-as)\,ete\,]+(e-as)(etf)+f\in eR\,e+eRf+fRf,
\end{eqnarray*}
which has a Peirce decomposition matrix of the form
$\footnotesize{\begin{pmatrix}k&q\\0&f\end{pmatrix}}$. Writing
down the Peirce matrix of $\,\alpha^{-1}\,$ in the form
$\footnotesize{\begin{pmatrix}w&x\\y&z\end{pmatrix}}$, we see that
$\,f=fz=z$, so $\,yq+zf=f\Rightarrow yq=0$. Since $\,yk=0\,$ and
$\,yf=0\,$ too, we have $\,y\alpha=0$, and so $\,y=0$. This
implies that $\,kw=e$. As $\,wk=e\,$ too, we conclude that
$\,k=a+(e-as)\,ete\in {\rm U}(eR\,e)$. This serves to show
that $\,{\rm sr}_{eRe}(a)=1$, as desired.\qed


\bs\nt
{\bf Remark 6.4.} In the setting of Theorem 6.2 above, the condition
$\,{\rm sr}_R(a)=1\,$ itself is in general not sufficient to imply
that $\,{\rm sr}_{eR\,e}(a)=1$.  For instance, if we take
$\,R={\mathbb M}_2\bigl({\mathbb Z}\bigr)\,$ with $\,e=E_{11}\,$
and $\,f=E_{22}$, then $\,a\!:=2E_{11}\in eR\,e\,$ has the property
that $\,{\rm sr}_R(a)=1\,$ by Theorem 5.8, but $\,{\rm sr}_{eR\,e}(a)
={\rm sr}_{\mathbb Z}(2)\neq 1\,$ as we have seen in Example 2.5B.


\bs
Using the the Suspension Theorem 6.2 in the special case where
$\,p=0$, we can easily retrieve the following well-known classical
corner ring result of Vaserstein [Va$_2$:~Theorem 2.8] for rings
of stable range one.

\bs\nt
{\bf Corollary 6.5.} {\it If $\,R\,$ is a ring of stable
range one, so is every corner ring $\,eR\,e\,$ of $\,R$.}

\bs\nt
{\bf Proof.} Let $\,f=1-e\in {\rm idem}\,(R)$, and consider any
$\,a\in eR\,e$.  Since $\,R\,$ has stable range one, we have
$\,{\rm sr}_R(a+f)=1$. By the Suspension Theorem, we will
automatically have $\,{\rm sr}_{eR\,e}(a)=1$. This shows that
the corner ring $\,eR\,e\,$ also has stable range one.\qed

\bs
In some sense, Corollary 6.5 is a ``companion'' to the result of Ara
and Goodearl in [AG], which implies that, in the case where $\,e\,$
is a {\it full idempotent\/} of $\,R\,$ (defined by the properties
that $\,e^2=e\,$ and $\,R\,eR=R\,$), the corner ring $\,eRe\,$ having
stable range one would guarantee that $\,R\,$ itself has stable
range one.  Combining this with Corollary 6.5, we see, for instance,
that a matrix ring $\,R={\mathbb M}_n(S)\,$ (for a fixed integer
$\,n$\/) has stable range one iff the base ring $\,S\,$ has stable
range one.  This is another well known result of Vaserstein in
[Va$_1$], which was also reiterated later in [Va$_2$:~Theorem 2.4].

\bs
By suitably applying the Suspension Theorem 6.2 and the
Product Theorem 2.8, we can now easily prove the following
``self-extensions'' of Theorem 6.2.

\bs\nt
{\bf Theorem 6.6.} {\it Let $\,e,f\in {\rm idem}\,(R)\,$ with
$\,e+f=1$.}

\mk\nt
(1) {\it If $\,a\in eR\,e$, $\,p\in fR\,e\,$ and $\,b\in fRf\,$
are such that $\,{\rm sr}_{eR\,e}(a)=1\,$ and $\,{\rm sr}_{fRf}(b)=1$,
then $\,{\rm sr}_R(a+p+b)=1$.}

\sk\nt
(2) {\it If $\,a\in eR\,e$, $\,p\in fR\,e$ and $\,u\in {\rm U}(fRf)$,
then $\,{\rm sr}_{eR\,e}(a)=1\,$ iff $\,{\rm sr}_R(a+p+u)=1$.}

\bs\nt
{\bf Proof.} (1) Since $\,{\rm sr}_{eRe}(a)=1$, the Suspension
Theorem 6.2 shows that $\,{\rm sr}_R(a+p+f)=1$.  Similarly,
$\,{\rm sr}_{fRf}(b)=1\,$ shows that $\,{\rm sr}_{R}(e+b)=1$. Since
$$
(a+p+f)\,(e+b)=ae+p\,e+fb=a+p+b,
$$
the Product Theorem 2.8 yields the desired conclusion that
$\,{\rm sr}_R(a+p+b)=1$.

\mk\nt
(2) The ``only if'' part follows from (1) and the trivial fact
that $\,{\rm sr}_{fRf}(u)=1$. For the ``if'' part, assume that
$\,{\rm sr}_R(a+p+u)=1$. Letting $\,v\,$ be the inverse of
$\,u\,$ in $\,fRf$, we have $\,{\rm sr}_R(e+v)=1\,$ since
$\,e+v\in {\rm U}(R)$.  By the Product Theorem 2.8, we have
$$
1={\rm sr}_R \bigl((a+p+u)\,(e+v)\bigr)={\rm sr}_R(a+p+f),
$$
so the Suspension Theorem 6.2 gives $\,{\rm sr}_{eR\,e}(a)=1$.\qed

\bs
Upon applying Theorem 6.6 to unit-regular elements (instead of just
stable range one elements), we obtain quickly the following result
on $\,{\rm ureg}\,(R)\,$ (for any ring $\,R\,$) which recovers
a part of the main theorem in the paper of Lam and Murray [LM:~\S2].

\bs\nt
{\bf Theorem 6.7.} {\it Let $\,e,f\in {\rm idem}\,(R)\,$ with
$\,e+f=1$, and let $\,a\in eR\,e$, $\,u\in {\rm U}(fRf)$. Then
$\,a+u\in {\rm ureg}\,(R)\,$ iff $\,a\in {\rm ureg}\,(eR\,e)$.}

\bs\nt
{\bf Proof.} The ``if'' part is trivial, so we need only prove the
``only if'' part.  Assuming that $\,a+u\in {\rm ureg}\,(R)$, we have
$\,{\rm sr}_{R}(a+u)=1\,$ by Theorem 4.1. Therefore, Theorem 6.6(2)
(applied with $\,p=0\,$) shows that $\,{\rm sr}_{eR\,e}(a)=1$. On the
other hand, if $\,r\,$ is any inner inverse of $\,a+u\,$ in $\,R$, we
see easily (by using Peirce's decomposition theorem) that $\,ere\,$ is
an inner inverse of $\,a\,$ in $\,eR\,e$, so $\,a\in {\rm reg}\,(eR\,e)$.
Coupling this with the information that $\,{\rm sr}_{eR\,e}(a)=1$, another
application of Theorem 4.1 shows that $\,a\in {\rm ureg}\,(eR\,e)$.\qed

\bs
Next, we will prove the theorem below on the stable range of certain
lower triangular matrices in matrix rings. In the special case where
such matrices have a zero diagonal, for instance, this result leads
to a large supply of nilpotent matrices with stable range one.

\bs\nt
{\bf Theorem 6.8.} {\it Let $\,R={\mathbb M}_n(S)$, where
$\,n\geq 1\,$ and $\,S\,$ is any ring. If $\,A=(a_{ij})\in R\,$
is a lower triangular matrix with $\,{\rm sr}_S(a_{ii})=1\,$
for all $\,i\in [1,\,n]$, then $\,{\rm sr}_R(A)=1$.}

\bs\nt
{\bf Proof.} The proof of this theorem will be carried out by
induction on $\,n\geq 1$. For $\,n=1$, the result is trivial.
For the inductive step, write $\,A\,$ in the block form
$\footnotesize{\begin{pmatrix} A'&0\\ \alpha & a_{nn}\end{pmatrix}}$,
where $\,A'\in {\mathbb M}_{n-1}(S)\,$ satisfies the same
hypothesis as $\,A$, and $\,\alpha\,$ is a row vector of length
$\,n-1$. By the inductive hypothesis, $\,A'\,$ has stable
range one in $\,{\mathbb M}_{n-1}(S)$. Applying Theorem 6.6(1)
(with $\,b=a_{nn}\in fRf\,$ having stable range one) for the
two complementary idempotents $\,e=E_{11}+\cdots +E_{n-1,n-1}\,$
and $\,f=E_{n,n}$, we see that $\,{\rm sr}_R(A)=1$. \qed

\bs\nt
{\bf Remark 6.9.} Needless to say, the theorem above also holds
for {\it upper\/} triangular matrices $\,A=(a_{ij})\,$ satisfying
the same hypothesis that $\,{\rm sr}_S(a_{ii})=1\,$ for all $\,i$.
After taking such a matrix $\,A\,$ and applying some row and
column permutations to it, we can come up with some examples
of possibly non-nilpotent matrices $\,A\in R={\mathbb M}_n(S)\,$
with $\,{\rm sr}_R(A)=1$. For instance, for $\,n=3$, we may
start with $\,A=\footnotesize{\begin{pmatrix}a&x&y\\0&b&z\\0&0&c
\end{pmatrix}}\in {\mathbb M}_3(S)$, where $\,x,y,z\in S\,$ are
arbitrary and (say) $\,\{a,b,c\}\subseteq {\rm idem}\,(S)\cup
{\rm U}(S)$.  Interchanging the first and third rows of
$\,A\,$ and then interchanging the first and second columns
of the resulting matrix, we see that any matrix of the
form $\footnotesize{\begin{pmatrix}0&0&c\\b&0&z\\x&a&y
\end{pmatrix}}$ has stable range one in $\,{\mathbb M}_3(S)\,$
for any base ring $\,S$.

\mk
For our next result, we will prove Theorem 6.10 below which shows
that what we have called ``Super Jacobson's Lemma 3.2'' in \S3 holds
not only for units, regular elements and unit-regular elements,
{\it but also for the set of stable range one elements\/} in an
arbitrary ring $\,S$.  The reason for presenting this result here
instead of in \S3 is that the proof of this new result for stable
range one elements in any ring $\,S\,$ makes substantial use of
the Suspension Theorem 6.2 in this section.

\bs\nt
{\bf Theorem 6.10.} {\it For any three elements $\,a,b,x \in S$,
$\,{\rm sr}_S(a+b-axb)=1\,$ if and only if $\,{\rm sr}_S(a+b-bxa)=1$.}
({\it In particular, after replacing $\,a,\,b\,$ by $\,1-a,\;1-b\,$
and letting $\,x=1$, it follows that $\,{\rm sr}_S(1-ab)=1\,$ iff
$\,{\rm sr}_S(1-ba)=1$.})

\bs\nt
{\bf Proof.} Working in the matrix ring $\,R={\mathbb M}_2(S)$,
we first observe that, by an elementary column operation
followed by an elementary row operation, we can bring the matrix
$\,M\!:=\footnotesize{\begin{pmatrix}1&a\\b&c\end{pmatrix}}$ first
to $\,\footnotesize{\begin{pmatrix}1&0\\b&c-ba\end{pmatrix}}$, and
then to $\,\footnotesize{\begin{pmatrix}1&0\\0&c-ba\end{pmatrix}}$.
Thus, the Suspension Theorem 6.2 shows that $\,{\rm sr}_R(M)=1\,$
iff $\,{\rm sr}_S(c-ba)=1$.  To prove Theorem 6.10, we note further
that, by an elementary column operation followed by a switch of rows,
we can bring a matrix $\,P=\footnotesize{\begin{pmatrix}1&-b\\1-ax&a
\end{pmatrix}}$ to $\,Q=\footnotesize{\begin{pmatrix}1&a\\1-bx&-b
\end{pmatrix}}$. Thus, $\,{\rm sr}_R(P)=1\,$ iff $\,{\rm sr}_R(Q)=1$.
Applying now the observation in the first part of our proof, we see
that $\,{\rm sr}_S(a-(1-ax)(-b))=1\,$ iff $\,{\rm sr}_S(-b-(1-bx)\,a)=1$.
After a sign change, we conclude that $\,{\rm sr}_S(a+b-axb)=1\,$
iff $\,{\rm sr}_S(a+b-bxa)=1$. \qed

\bs\nt
{\bf Example 6.11.} In contrast to the last statement in
Theorem 6.10, it is not difficult to produce a ring $\,S\,$
with three elements $\,a,b,c\,$ such that $\,{\rm sr}_S(c-ab)=1\,$
but $\,{\rm sr}_S(c-ba)\neq 1$.  To get such an example, we take
$\,S={\mathbb M}_2({\mathbb Z})$, $\,a=E_{12}$, $\,b=E_{11}$ and
$\,c=2E_{21}\,$ in $\,S$.  With these choices, $\,c-ab=c=2E_{21}\,$
has stable range one by Theorem 5.8. On the other hand,
$\,c-ba=\footnotesize{\begin{pmatrix}0&-1\\2&0\end{pmatrix}}$
is matrix-equivalent to $\,{\rm diag}\,(2,1)$, which
according to Theorem 6.6(2) does not have stable range one
in $\,S={\mathbb M}_2({\mathbb Z})\,$ since
$\,{\rm sr}_{\mathbb Z}(2)\neq 1$.

\bs
To give an application of Example 6.11, we trace our path back to
Corollary 3.10, which proved the invariance of the stable range one
property with respect to the transposition of $\,n\times n\,$
matrices {\it over a commutative ring $\,S$.}  In the following,
we will offer an example to show that the commutativity assumption
on $\,S\,$ cannot be dropped, even in the case $\,n=2$.  This
should not be a huge surprise if we recall, for instance from
[Ja$_1$: p.\,24] (c.~1953), that the transpose of a $\,2\times 2\,$
invertible matrix over a noncommutative ring need not be invertible
in general.  Such a classical fact suggests that we should not
expect a whole lot about the preservation of matrix properties
under transposition when we work with matrices over noncommutative
rings, as the following concrete example shows.

\bs\nt
{\bf Example 6.12.} Let $\,S\,$ be a ring with three elements $\,a,b,c\,$
such that $\,{\rm sr}_{S}(c-ab)=1$, but $\,{\rm sr}_{S}(c-ba)\neq 1$.
(For instance, we can take the ring $\,S\,$ that was constructed in
Example 6.11.)  In the matrix ring $\,R\!:={\mathbb M}_2(S)$, let
$\,M\!:=\footnotesize{\begin{pmatrix}1&a\\b&c \end{pmatrix}}$, with
transpose $\,M^T=\footnotesize{\begin{pmatrix}1&b\\a&c \end{pmatrix}}$.
Using the criterion for $\,{\rm sr}_R(M)=1\,$ obtained in the proof
of Theorem 6.10 (and its analogue for $\,M^T\,$), we see that
$\,{\rm sr}_R(M)=1$, but $\,{\rm sr}_R(M^T)\neq 1$.

\bs\nt
{\bf Example 6.13.} Let $\,R\!:={\mathbb M}_4(\mathbb{Z})$ and $\,M\!:=\footnotesize{\begin{pmatrix}1&a\\b&c \end{pmatrix}} \in R$ with $a,\,b,\,c$ as in Example 6.11. We saw in Example 6.12 that $\,{\rm sr}_R(M)=1\,$. Note that $\,{\rm det}\,(M)=2$  and $\,{\rm sr}_{\mathbb Z}(2)\neq 1\,$. This shows that the determinant of a matrix with stable range one over a commutative ring may not have stable range one, in answer to a question raised by a referee of this
paper.

\mk

\bs\nt
{\bf \S7. Integral Matrices of Stable Range One}

\bs
By ``integral matrices'' in the above caption, we mean the matrices
in the rings $\,{\mathbb M}_n\bigl({\mathbb Z}\bigr)$. In this
section, we shall determine precisely the set of all such matrices
with stable range one, for any $\,n\geq 1$.  To begin with, we
first study the case of {\it diagonal\/} matrices $\,A\,$ over any
ring $\,S$.  In the case where some diagonal entry of $\,A\,$ is
zero, we have the following slightly surprising result.


\bs\nt
{\bf Theorem 7.1.} {\it Let $\,A={\rm diag}\,(a_1,\dots,a_n)\in
R={\mathbb M}_n(S)$, with some $\,a_i=0$. Then $\,{\rm sr}_R(A)=1$.}

\bs\nt
{\bf Proof.} We may assume, of course, that $\,n\geq 2$. Without loss
of generality, we may also assume that $\,a_n=0$. After switching the
first and the last rows, we obtain a matrix $\,B\,$ with last row
$\,(a_1,0,\dots,0)$, and with northwest $\,(n-1)\times (n-1)\,$ block
given by $\,D={\rm diag}\,(0,a_2,\dots,a_{n-1})$. By invoking an
inductive hypothesis here, we may assume that $\,{\rm sr}\,(D)=1\,$
in $\,{\mathbb M}_{n-1}(S)$.  Applying now Theorem 6.6(1) (with
respect to the idempotent $\,e=E_{11}+\cdots+E_{n-1,n-1}$), we see
that $\,{\rm sr}_R(A)={\rm sr}_R(D)=1$, as desired. \qed

\bs
In the paper [CP$_1$] of C\v{a}lug\v{a}reanu and Pop, $2\times 2$ matrices
over $\,{\mathbb Z}\,$ of stable range one were completely determined;
namely, these were shown to be precisely the $\,2\times 2\,$ integral
matrices $\,A\,$ with $\,{\rm det}\,(A)\in \{0\}\cup {\rm U}({\mathbb Z})
=\{0,\,{\pm 1}\}$.  With the aid of Theorem 7.1, we can now extend
this theorem to all $\,n\times n\,$ matrices over $\,{\mathbb Z}$.
To do this, recall that according to Henriksen [He$_1$], a ring $\,S\,$
is called an {\bf elementary divisor ring\/} if every square matrix
over $\,S\,$ is equivalent to a diagonal matrix; that is, for any
$\,A\in {\mathbb M}_n(S)$, there exist $\,U,V\in {\rm GL}_n(S)\,$
such that $\,UAV\,$ is a diagonal matrix. For instance, it is well
known that every commutative PID (principal ideal domain) is an
elementary divisor ring.

\bs\nt
{\bf Theorem 7.2.} {\it Let $\,A\in R={\mathbb M}_n(S)$, where
$\,n\geq 1\,$ and $\,S\,$ is a commutative elementary divisor
domain. If $\,{\rm det}\,(A)\in \{0\}\cup {\rm U}(S)$, then
$\,{\rm sr}_R(A)=1$.\footnote{This implies, in particular, that
any nilpotent matrix in $\,R\,$ has stable range one.} In the
special case where $\,S={\mathbb Z}$, the converse of this
statement also holds.}

\bs\nt
{\bf Proof.} To prove the first conclusion, we may assume (in view
of the Product Theorem 2.8) that $\,A={\rm diag}\,(a_1,\dots,a_n)$.
If $\,{\rm det}\,(A)\in \{0\}\cup {\rm U}(S)$, then either some
$\,a_i=0\,$ or $\,A\in {\rm U}(R)$.  In the latter case, we have
trivially $\,{\rm sr}\,(A)=1$.  In the former case, the same
conclusion holds also in view of Theorem 7.1.  To prove the last
statement in the theorem, we work now in the case $\,S={\mathbb Z}$,
and assume that $\,{\rm sr}\,(A)=1$. Replacing $\,a_i\,$ by $\,-a_i\,$
if necessary, we may assume that $\,a_i\geq 0\,$ for all $\,i$.
If $\,d\!:={\rm det}\,(A)=a_1\cdots a_n\notin \{0\}\cup {\rm U}(S)
=\{0,\,{\pm 1}\}$, then $\,d\geq 2$. Letting
$$
D={\rm diag}\,(a_1,\dots,a_n)\cdot {\rm diag}\,(a_2,\dots,a_n,a_1)
\cdots {\rm diag}\,(a_n,a_1,\,\dots\,a_{n-1})=d\cdot I_n\,,\leqno (7.3)
$$
we have $\,{\rm sr}\,(D)=1\,$ by Theorem 2.8, since each factor of
$\,D\,$ above has stable range one. Applying the stable range one
criterion in Theorem 2.2(3), we see that there exists a matrix
$\,B\in {\mathbb M}_n\bigl({\mathbb Z}\bigr)\,$ such that
$$
M\!:= D+\bigl(I_n+D\cdot d^n I_n\bigr)\cdot B
= D+(1+d^{n+1})\,B \in {\rm GL}_n\bigl({\mathbb Z}\bigr).
\leqno (7.4)
$$
Computing $\,{\rm det}\,(M)\,$ over the factor ring
$\,{\mathbb Z}/(1+d^{n+1})\,{\mathbb Z}\,$ and noting that
$\,{\rm det}\,(M)={\pm 1}$, we get a congruence relation
$\,{\rm det}\,(D)=d^n \equiv {\pm 1}\;
\bigl({\rm mod}\;(1+d^{n+1})\bigr)$. Since $\,d\geq 2$,
this is clearly impossible, so we have proved that
$\,{\rm det}\,(A)\in \{0,\,{\pm 1}\}$, as desired. \qed

\bs
By using the theorem above in the special case $\,S={\mathbb Z}$,
any issue concerning integral matrices of stable range one can
be easily settled through a simple computation of determinants.
The following consequence of Theorem 7.2 gives some good
instances of this.

\bs\nt
{\bf Corollary 7.5.} {\it For the matrix ring
$\,R={\mathbb M}_n\bigl({\mathbb Z}\bigr)\,$ where $\,n\geq 1$,
the following holds.}

\sk\nt
(1) {\it A diagonal matrix $\,{\rm diag}\,(a_1,\dots,a_n)\in R\,$
has stable range one iff $\,a_i=0\,$ for some $\,i$, or
$\,a_i\in \{{\pm 1}\}\,$ for all $\,i$.}\\
(2) {\it For any two matrices $\,A,B\in R$, $\,{\rm sr}_R(AB)=1\,$
iff $\,{\rm sr}_R(BA)=1\,$.}


\bs\nt
{\bf Proof.} (1) follows directly from Theorem 7.2, while
(2) also follows from Theorem 7.2 since
$\,{\rm det}\,(AB)={\rm det}\,(A)\cdot {\rm det}\,(B)={\rm det}\,(BA)$. \qed


\bs
Coming back once more to the last ``converse'' statement
in Theorem 7.2, we should note that, while that statement
holds in the special case $\,S={\mathbb Z}$, it may not hold if
$\,S\,$ is just a commutative PID. For instance, if $\,S\,$ is a
commutative semilocal PID (e.g.~a discrete valuation domain), it
is well known that both $\,S\,$ and $\,R={\mathbb M}_n(S)\,$ are
rings of stable range one.  In this case, $\,{\rm sr}_R(A)=1\,$
holds {\it for all $\,A\in R$,} but we certainly may not have any
conclusion such as $\,{\rm det}\,(A)\in \{0\}\cup {\rm U}(S)$.
Finally, we should also point out that, if $\,S\,$ is commutative
domain and $\,A\in R={\mathbb M}_n(S)\,$ has $\,{\rm det}\,(A)=0$,
we may not have $\,{\rm sr}_R(A)=1\,$ either.  For instance, if
$\,S\,$ is a polynomial ring $\,k\,[x,y]\,$ over a field $\,k$,
the matrix $\,A=\footnotesize{\begin{pmatrix}x&y\\0&0\end{pmatrix}}$
with $\,{\rm det}\,(A)=0\,$ will not have stable stable range one
in $\,{\mathbb M}_2(S)$.  The theoretical justification for this
statement will be more fully discussed in our later work.

\mk

\bs\nt
{\bf \S8. \ Open Questions}

\bs
In spite of the work done in this paper, there are still many
significant unsolved problems in the study of stable range one
elements in rings. To stimulate further work in this area, we
collect and formulate three such open problems below.

\bs
The first open problem is prompted by the fact that we have not
been able to fully understand the relationship between nilpotent
elements, strongly $\pi$-regular elements, and elements of stable
range one in an arbitrary ring $\,R$, in spite of some of our
partial results in \S5 and \S6.

\bs\nt
{\bf Question 8.1.} {\it For any given nilpotent or strongly
$\,\pi$-regular element $\,a\in R$, what is a necessary and
sufficient condition for $\,{\rm sr}_R(a)=1\,$}?

\bs
From Theorem 5.1, Theorem 5.2 and Theorem 6.8, we do have some
partial answers to the question above if we impose some further
assumptions on the element $\,a\in R\,$ or on the ring $\,R$.
A general answer for Question 8.1, even just for an exchange
ring $\,R$, should lead to an interesting overview for some of
these earlier results in \S5 and \S6. In the meantime, any such
characterization results in the case of strongly $\pi$-regular
elements should be helpful toward a further understanding of
Ara's highly significant result [Ar:~Theorem 4] that all
strongly $\,\pi$-regular rings have stable range one.

\bs
A second question we want to raise concerns the behavior of stable
range one elements in a ring $\,R\,$ with respect to the passage
from $\,R\,$ to its various factor rings.  In their 1974 paper [FS],
Fisher and Snider proved essentially that a ring element $\,a\in R\,$
is strongly $\,\pi$-regular iff its images $\,\overline{a}=a+P\,$
are strongly $\,\pi$-regular in $\,R/P\,$ for every prime ideal
$\,P\subset R$.  In view of this classical result, it would be
tempting to ask the following.

\bs\nt
{\bf Question 8.2.} {\it If a ring element $\,a\in R\,$ is such that
its images $\,\overline{a}=a+P\,$ have stable range one in $\,R/P\,$
for every prime ideal $\,P\subset R$, does it follow that
$\,{\rm sr}_R(a)=1\,$}?

\bs
Our third question concerns the general behavior of the
element-wise stable range one notion under ring extensions.

\bs\nt
{\bf Question 8.3.} {\it Let $\,a\,$ be a given element in
a ring $\,R$.}

\mk\nt
(A) {\it When does there exist a ring $\,K\supseteq R\,$\/}
({\it with $\,1_K=1_R\,$}) {\it such that $\,{\rm sr}_K(a)=1\,$}?

\sk\nt
(B) {\it When is it true that $\,{\rm sr}_K(a)=1\,$ for every
ring $\,K\supseteq R\,$\/} ({\it with $\,1_K=1_R\,$})?

\bs
Classically, if $\,R\,$ is a commutative integral domain, the answer
to Question (A) is trivially ``always'', as we can take $\,K\,$
to be the quotient field of $\,R$. More generally, if $\,R\,$ is
a {\it commutative\/} ring, the answer to Question (A) is still
``always'' since a result of Goodearl and Menal [GM:~Prop.~7.6]
guarantees that $\,R\,$ can be unitally embedded into a commutative
ring $\,K\,$ of stable range one. If $\,R\,$ is a {\it noncommutative\/}
ring, however, the situation would become more complicated. In the
case where $\,a\,$ is a product of units and idempotents in $\,R\,$
(in any order), we will surely have $\,{\rm sr}_K(a)=1\,$ in every
ring $\,K\supseteq R\,$ (with $\,1_K=1_R\,$).  On the other hand,
suppose $\,a,\,b\in R\,$ are such that $\,ab=1\neq ba$.  Then for
every ring $\,K\supseteq R\,$ (with $\,1_K=1_R\,$), we cannot
have $\,{\rm sr}_K(a)=1$.  Indeed, if $\,{\rm sr}_K(a)=1$, then
$\,a=aba \in {\rm reg}\,(K)\Rightarrow a\in {\rm ureg}\,(K)\,$
according to Theorem 4.1.  This would imply that $\,a\in {\rm U}(K)$,
and hence $\,ab=1\in K \Rightarrow ba=1\in K$, a contradiction.



\sk

\bs\nt
{\bf Acknowledgements.} It is the authors' great pleasure to thank Professors Pere Ara, Ebrahim Ghashghaei, Kenneth Goodearl, Pace Nielsen, and two anonymous referees, for a number of insightful comments and suggestions on this paper.

\sk

\bs\nt

\bs
\nt Department of Mathematics \\
\nt Panjab University \\
\nt Chandigarh-160014, India

\sk
\nt {\tt dkhurana@pu.ac.in}

\bs
\nt Department of Mathematics \\
\nt University of California \\
\nt Berkeley, CA 94720, USA

\sk
\nt {\tt lam@math.berkeley.edu}

\end{document}